\documentclass{amsart}

\usepackage{diagrams}
\usepackage[dvips]{graphicx}
\usepackage{epic}
\usepackage{latexsym}
\usepackage[T1]{fontenc}
\usepackage{amsrefs}
\usepackage{amssymb,amsmath}
\usepackage{amsthm}
\usepackage{mathrsfs}
\usepackage{enumerate}
\usepackage[all]{xy}

\newcommand{\id}{1\!\!1}
\newcommand{\geqs}{\geqslant}
\newcommand{\leqs}{\leqslant}
\newcommand{\mquad}{\!\!\!\!\!\!}

\newcommand{\lra}{\longrightarrow}

\newcommand{\rar}{\rightarrow}

\newarrow{TTo} ----{->}
\newarrow{Mapstto} |---{->}
\newarrow{Intto} C---{->}

\newtheorem{thm}{Theorem}[section]
\newtheorem*{thm*}{Theorem}
\newtheorem{defi}[thm]{Definition}
\newtheorem{lem}[thm]{Lemma}
\newtheorem{prop}[thm]{Proposition}	
\newtheorem*{prop*}{Proposition}	
\newtheorem{coro}[thm]{Corollary}
\newtheorem{exm}[thm]{Example}
\newtheorem{obs}[thm]{Observation}
\newtheorem*{obs*}{Observation}

\begin{document}

\title{OCHA and the swiss-cheese operad}        % Title of Document

\author{Eduardo Hoefel}                         % First Author
\email{hoefel@ufpr.br}                          % First Author's email
\address{Dep. de Matemática\\                   % First Author's postal address
         Universidade Federal do Paran\'a\\
         Brasil\\
         Curitiba\\
         c.p. 019081  cep: 81531-990}

%\keywords{Coloured operads, Axelrod-Singer compactification,
%homotopical algebra, resolutions of operads.}

%\classification{18G55, 18D50.}

\thanks{The author wishes to thank Jim Stasheff and Murray Gerstenhaber for the kind hospitality
during his stay as a visiting graduate student at the University of Pennsylvania (CNPq-Brasil grant
SWE-201064/04). We are also grateful to J. Stasheff and H. Kajiura for valuable discussions.
Some results of this paper consist of strengthened versions of some of the results present in the author's
Ph.D. thesis \cite{Hoe06c} defended in 2006 at Unicamp under the supervision of A. Rigas
and T. E. Barros.}

\begin{abstract}
In this paper we show that the relation between Kajiura-Stasheff's OCHA and
A. Voronov's swiss-cheese operad is analogous to the relation between
SH Lie algebras and the little discs operad. More precisely, we show that the OCHA
operad is quasi-isomorphic to the operad generated by the top-dimensional homology classes
of the swiss-cheese operad.
\end{abstract}

%\received{February 05, 2008}
%\revised{January 01, 2009}
%\published{April 28, 2009}
%\submitted{James Stasheff}
%%%arXiv: 0710.3546
%\volumeyear{2009}
%\volumenumber{4}
%\issuenumber{1}
%
%\startpage{83}

\maketitle

\section*{Introduction}
OCHA refers to the homotopy algebra of open and closed strings introduced by Kajiura and Stasheff \cite{KS06a}
inspired by the work of Zwiebach on string field theory \cite{Zwi98}. In \cite{KS06a} the
{\it $A_\infty$-algebras over $L_\infty$-algebras} are also introduced, they are the strong homotopy version
of $\mathfrak{g}$-algebras (or Leibniz pairs, see \cite{FGV95}). An OCHA is a structure obtained
by adding other operations to an $A_\infty$-algebra over an $L_\infty$-algebra. The physical meaning
of those additional operations is given by the ``opening of a closed string into an open one''.

Considering that its relevance to Physics is well acknowledged (see also \cite{KS06b}), in the present paper
we further explore the {\it mathematical significance} of a full OCHA, not restricted to an
$A_\infty$-algebra over an $L_\infty$-algebra. In \cite{Hoe06b} we have proven that any degree one coderivation
$D \in {\rm Coder}(S^c L \otimes T^cA)$ such that $D^2 = 0$ defines an OCHA structure on the pair
$(L,A)$. In this work we study the relation between OCHA's and A. Voronov's swiss-cheese operad and show
that it is analogous to the relation between SH Lie algebras and the little discs operad.
A graded Lie algebra is part of the structure of a Gerstenhaber Algebra,
which in turn is equivalent to an algebra over the homology little discs operad.
The Lie part of a Gerstenhaber algebra is given by the top-dimensional homology classes
of the little discs operad.

We will study the suboperad of the homology swiss-cheese operad generated by top-dimensional
homology classes and show that it is quasi-isomorphic to the operad whose algebras are OCHAs.
The quasi-isomorphism, however, is not of operads but only a quasi-isomorphism
of modules over the operad $\mathcal{L}_\infty$ of $L_\infty$-algebras.

Let $\mathcal{OC}_\infty$ be the OCHA operad and let $\mathcal{OC}$ denote the suboperad
of the homology swiss-cheese operad generated by top-dimensional homology classes.
Our main result is the following.
\begin{thm*}
There is a morphism of differential graded $\mathcal{L}_\infty$-modules
$\mu : \mathcal{OC}_\infty \lra \mathcal{OC}$
which induces an isomorphism in cohomology.
\end{thm*}
The paper is organized as follows.
In section \ref{homology_ld} we briefly review F. Cohen's theorem on the homology of the little discs
operad and state it using trees. Section \ref{homology_sc} reviews the analogous description of
the homology swiss-cheese operad in terms of generators and relations given by trees. In section
\ref{OCHA} we define OCHA in a grading and signs convention which is different from the original
one in \cite{KS06a}. The definition given here is appropriate for studying its correspondence with
the compactified configuration space of points on the closed upper half plane. We show that both definitions
are equivalent through the (de)suspension operator. A definition of the OCHA operad $\mathcal{OC}_\infty$ is
provided in section \ref{OCHA_operad} using the partially planar trees, a type of tree which is defined in the
same section. Section \ref{c(p,q)} reviews the construction of the compactified configuration spaces
$\overline{C(p,q)}$ first introduced by Kontsevich in \cite{Kont03}. The combinatorial structure of its
boundary strata is described in terms of partially planar trees and some examples are provided.
The well known equivalence between $\overline{C(1,q)}$ and the cyclohedron $W_{q+1}$
is explained in terms of those trees. In section \ref{OCHA_SS} we prove the quasi-isomorphism between $\mathcal{OC}_\infty$ and $\mathcal{OC}$ viewed as modules over the operad $\mathcal{L}_\infty$ of
$L_\infty$ algebras.
The main tool used in its proof is the spectral sequence of $\overline{C(p,q)}$ as a manifold with corners.

\subsection*{Notation and Conventions}
Let us fix a field $k$ of characteristic zero. In this paper, all vector spaces are over $k$ and
`graded vector space' will always mean `$\mathbb{Z}$-graded vector space', unless otherwise stated.
Let $V$ be a graded vector space, we define a left action of the symmetric group $S_n$ on $V^{\otimes n}$ in the
following way:  if $\tau \in S_2$ is a transposition, then the action is given by
$x_1 \otimes x_2 \stackrel{\tau}{\mapsto} (-1)^{|x_1||x_2|} x_2 \otimes x_1$. Since any $\sigma \in S_n$ is a
composition of transpositions, the sign of the action of $\sigma$ on $V^{\otimes n}$ is well defined:
\begin{equation}
 x_1 \otimes \cdots \otimes x_n \stackrel{\sigma}{\mapsto} \epsilon (\sigma)
x_{\sigma(1)} \otimes \cdots \otimes x_{\sigma(n)}.
\end{equation}
We will refer to $\epsilon (\sigma)$ as the Koszul sign of the permutation. Let us
define $\chi(\sigma) = (-1)^\sigma \epsilon (\sigma)$, where $(-1)^\sigma$ is the
sign of the permutation.

Given two homogeneous maps $f,g : V \rar W$ between graded vector spaces,
according to the Koszul sign convention (which will be used throughout this work), we have:
\begin{equation} \label{koszul}
(f \otimes g)(v_1 \otimes v_2) = (-1)^{|g||v_1|} (f(v_1) \otimes g(v_2)).
\end{equation}

We will use the notation of Lada-Markl \cite{LM95}
for the suspension and desuspension operators:
$\uparrow$ and $\downarrow$. Let
$\uparrow\!\! V$ (resp. $\downarrow\!\! V$) denote the suspension (resp. desuspension)
of the graded vector space $V$ defined by:
$(\uparrow\!\! V)^p = V^{p-1}$ (resp. $(\downarrow\!\! V)^p = V^{p+1}$).  We thus have the natural
maps $\uparrow : V \rar \uparrow\!\! V$ of degree $1$, and $\downarrow : V \rar \downarrow\!\! V$ of degree $-1$.
Let $\uparrow^{\otimes n}$ denote $\bigotimes^n \uparrow : \bigotimes^n V \rar \bigotimes^n \uparrow\!\! V$
and $\downarrow^{\otimes n}$ is defined analogously.
The operators $\uparrow^{\otimes n}$ and $\downarrow^{\otimes n}$ transform symmetric operations into anti-symmetric ones. In fact, let $E$ (resp. $A$) denote the symmetric
(resp. anti-symmetric) left action of the group of permutations $S_n$ on $V^{\otimes n}$:
\begin{align}
E(\sigma)(x_1 \otimes \cdots \otimes x_n)& =
\epsilon (\sigma) x_{\sigma(1)} \otimes \cdots \otimes x_{\sigma(n)} \label{E}\\
A(\sigma)(x_1 \otimes \cdots \otimes x_n)& =
\chi (\sigma) x_{\sigma(1)} \otimes \cdots \otimes x_{\sigma(n)} \label{A}
\end{align}

Both actions are related by:
$\uparrow^{\otimes n} E(\sigma) \downarrow^{\otimes n} = (-1)^{n(n-1)/2} A(\sigma)$, for any
$\sigma \in S_n$. In particular, $\uparrow^{\otimes n} \circ \downarrow^{\otimes n} = (-1)^{n(n-1)/2} \cdot \id$. The sign $(-1)^{n(n-1)/2}$ is a consequence of the Koszul sign convention (\ref{koszul}) defined above
(see also \cite{DMZ07}).

Let us now describe how the notation for operads and its related concepts (such as: representations, ideals and modules)
will be used in this paper. Our description will not necessarily include precise definitions. Those can be found
in \cite{MSS02}. An operad is any sequence $\mathcal{O} = \{ \mathcal{O}(n) \}_{n \geqs 1}$ of objects in a symmetric
monoidal category (such as the category of topological spaces or the category of vector spaces)
endowed with a right action of the symmetric
group $S_n$ on each $\mathcal{O}(n)$ and a composition law satisfying natural associative and equivariance conditions.

Given a graded vector space $V$, the endomorphism operad of $V$ is defined as
${\rm End}_{V}(n) = {\rm Hom}(V^{\otimes n},V)$. The composition law $\circ_i$ in ${\rm End}_{V}$ is
defined by the usual composition in the $i$th variable of multilinear maps and the right action of
$S_n$ on ${\rm Hom}(V^{\otimes n},V)$ is the composition with the symmetric left action $E$ defined
by (\ref{E}). In particular, this means that for graded vector spaces, according to our conventions, `symmetric'
always mean `graded symmetric'.

Among the standard examples of operads are those defined in terms of {\it trees}.
In this paper, in accordance with \cite{GJ94}, trees are oriented and not necessarily compact:
an edge may be terminated by a vertex at only one end (or none).
Such an edge is called {\it external}. An external edge oriented toward
its vertex is called a {\it leaf}, otherwise it is called {\it the root}. Trees
are assumed to have only one root. The leaves of each tree are labeled by natural numbers.
The action of $S_n$ on trees with $n$ leaves is defined by permuting the labels.
The composition law $\circ_i$ on operads defined in terms of trees is given by the grafting operation,
i.e., the identification of the root of one tree with the leaf labelled $i$ of the other tree.

We also need to mention the {\it coloured operads}, a concept that goes back to Boardman and Vogt
\cite{BV73}. Following the notation of Berger and Moerdijk \cite{BM06}, given a set of colours $C$, a $C$-coloured
operad $\mathcal{P}$ is defined by assigning to each $(n+1)$-tuple of colours $(c_1, \dots, c_n;c)$ an
object \[ \hskip 10em \mathcal{P}(c_1, \dots, c_n;c) \qquad \mbox{ in some monoidal category} \]
endowed with a composition law and a symmetric group action. The defining conditions for coloured operads
are analogous to those for ordinary operads.

Given a family  $A = \{ A_c \}_{c \in C}$ of vector spaces indexed by $C$, the $C$-coloured operad
${\rm End}(A)$ is defined by:
\begin{equation}\label{end_C}
 {\rm End}(A)(c_1, \dots, c_n; c) = {\rm Hom}(A_{c_1} \otimes \cdots \otimes A_{c_n},A_{c}).
\end{equation}
For coloured operads, the composition law is only defined when the colour of the output coincides with
the colour of the input. Another example of a coloured operad is given by trees with coloured edges, i.e.,
trees such that for each edge is assigned a element in some set $C$. For trees with coloured edges, the grafting
operation is only defined when the colour of the root coincides with the colour of the corresponding leaf.
Coloured trees will be used throughout the present paper. In fact, we will use 2-coloured trees where the
colours of the edges are wiggly or straight, according to the notation used in \cite{KS06a}.

Let $\mathcal{P}$ be an operad and let $M = \{ M(n) \}_{n \geqs 1}$ be a sequence of objects
where each $M(n)$ has a right $S_n$-action. We say that $M$ is a left $\mathcal{P}$-module if it is
endowed with a left `operadic action' $\circ_i^\lambda$:
\begin{equation}\label{left_module}
 \circ_i^\lambda: \mathcal{P}(n) \otimes M(m) \rar M(m + n - 1)
\end{equation}
which is equivariant and satisfies associative conditions analogous to those in the definition of operads.
The definition of right modules is similar.

An {\it ideal} in an operad $\mathcal{P}$ is a sequence of objects $\mathcal{I} = \{ I(n) \}_{n \geqs 1}$
with $I(n) \subseteq \mathcal{P}(n)$, where each $I(n)$ is invariant under the action of $S_n$
and $\mathcal{I}$ is a left and right $\mathcal{P}$-module. We refer the reader to \cite{MSS02} for
the precise definitions and further details about these concepts.

\section{The homology little disks operad}\label{homology_ld}
We begin by recalling the description of the homology little discs operad $H_\bullet (\mathcal{D})$
in terms of generators and relations.  %\subsection{Generators and relations for $H_\bullet (D)$}
To keep our notation in accordance with \cite{KS06a}, we
will represent classes in $H_\bullet (\mathcal{D})$ by trees with wiggly edges.
As usual in operad theory, all trees are assumed to be {\it rooted} and {\it oriented} toward the root.
A tree with only one vertex and $n$ incoming edges is called an $n$-corolla.

The little disks operad (also called little 2-disks operad) is a sequence
$\mathcal{D} = \{ \mathcal{D}(n) \}_{n \geqs 1}$ of topological spaces $\mathcal D(n)$ defined as the space
of all maps
\[ d : \coprod_{1 \leqs s \leqs n} D_s \rar D \]
from the disjoint union of $n$ numbered standard two dimensional disks $D_1, \dots, D_n$ to $D$, where
\[ D := \{ (x_1,x_2) \in \mathbb{R}^2 \ | \ x_1^2 + x_2^2 \leqs 1 \}, \]
such that $d$, when restricted to each disk, is a composition of translation and multiplication by a positive
real number and the images of the interiors of the disks are disjoint. The symmetric group acts by renumbering
the disks. We may interpret $d \in \mathcal{D}(n)$ as the standard disk $D$ with $n$ numbered disjoint circular
holes. The operad composition $\gamma_{\mathcal D}(d;d_1, \dots, d_n)$
is, intuitively speaking, given by gluing $n$ disks $d_1, \dots, d_n$
in the holes of the disk with holes $d$, and erasing the seams (see \cite{MSS02,May72} for more precise
definitions).

Since $\mathcal{D}(2)$ is homotopy equivalent
to $S^1$, its homology has two generators. The zero dimensional generator will be denoted
by a 2-corolla with a ``white'' vertex: $\circ$, while the one dimensional generator will be denoted by a 2-corolla with a ``black'' vertex: $\bullet$. Let
\[ \Bigg\{ \ \raisebox{-1.5em}{\begin{picture}(0,0)%
\includegraphics{sm2_tree.pstex}%
\end{picture}%
\setlength{\unitlength}{3947sp}%
\begingroup\makeatletter\ifx\SetFigFont\undefined%
\gdef\SetFigFont#1#2#3#4#5{%
  \reset@font\fontsize{#1}{#2pt}%
  \fontfamily{#3}\fontseries{#4}\fontshape{#5}%
  \selectfont}%
\fi\endgroup%
\begin{picture}(665,702)(636,-241)
\put(626,380){\makebox(0,0)[lb]{\smash{{\SetFigFont{8}{9.6}{\rmdefault}{\mddefault}{\updefault}$1$}}}}
\put(1107,380){\makebox(0,0)[lb]{\smash{{\SetFigFont{8}{9.6}{\rmdefault}{\mddefault}{\updefault}$2$}}}}
\end{picture}%
},
\raisebox{-1.5em}{\begin{picture}(0,0)%
\includegraphics{sl2_tree.pstex}%
\end{picture}%
\setlength{\unitlength}{3947sp}%
\begingroup\makeatletter\ifx\SetFigFont\undefined%
\gdef\SetFigFont#1#2#3#4#5{%
  \reset@font\fontsize{#1}{#2pt}%
  \fontfamily{#3}\fontseries{#4}\fontshape{#5}%
  \selectfont}%
\fi\endgroup%
\begin{picture}(665,702)(636,-241)
\put(626,380){\makebox(0,0)[lb]{\smash{{\SetFigFont{8}{9.6}{\rmdefault}{\mddefault}{\updefault}$1$}}}}
\put(1107,380){\makebox(0,0)[lb]{\smash{{\SetFigFont{8}{9.6}{\rmdefault}{\mddefault}{\updefault}$2$}}}}
\end{picture}%
} \Bigg\} \] be a basis for $H_\bullet (\mathcal{D}(2))$, where the
first generator has degree zero, the second has degree one and both are invariant under the action
of the permutation group. An algebra over $H_\bullet (\mathcal{D})$ will thus have two graded symmetric
operations. Notice however that the bracket defined below by (\ref{gbracket})
is {\it skew} graded symmetric, as required by the definition of a Gerstenhaber algebra.
F. Cohen's famous result about the homology of $\mathcal{D}$ can be stated
in the language of trees, as follows.
\begin{thm}[F. Cohen \cite{Cohen76}] \label{cohen76}
The homology little disks operad $H_\bullet (\mathcal{D})$ is isomorphic, as an operad of $\mathbb{Z}$-graded vector spaces, to the operad generated by the above trees, subject to the following relations:
\begin{enumerate}[i)]
\item Both generators are invariant under permutation of their labels;
\item Jacobi identity: \raisebox{-2em}{\input{jacobi_tree.pstex_t}};
\item Leibniz rule: \raisebox{-2em}{\input{leib_tree.pstex_t}}.
\end{enumerate}
\end{thm}
Let us now describe algebras over $H_\bullet (\mathcal{D})$. It is well known that algebras over
$H_\bullet (\mathcal{D})$ are equivalent to Gerstenhaber algebras. However, a brief description
of that equivalence will help clarify our exposition involving the swiss-cheese
operad and its relations to OCHA.

Remember that $V$ is an algebra
over $H_\bullet (\mathcal{D})$ if there is a morphism of operads $\Phi : H_\bullet (\mathcal{D}) \rar {\rm End}_{V}$ where
${\rm End}_{V}(n) = \{ {\rm Hom}(V^{\otimes n},V) \}$, is the endomorphism operad.
Consequently, $V$ is a graded vector space endowed with two graded symmetric
operations $m_2 : V \otimes V \rar V$ of degree $|m_2| = 0$, corresponding to the
first generating tree, and $l_2 : V \otimes V \rar V$ of degree $|l_2| = 1$
corresponding to the second generating tree.

The two relations presented above in terms of trees correspond to the equalities:
\begin{equation} \label{first}
\displaystyle{\sum_\sigma l_2 \circ (\id \otimes l_2) \circ E(\sigma) = 0}
\end{equation}
and
\begin{equation} \label{second}
l_2 \circ (\id \otimes m_2) = m_2 \circ (l_2 \otimes \id) +
                         m_2 \circ (\id \otimes l_2) \circ E(\tau_{1,2}).
\end{equation}
In the first equality, $\sigma$ runs over all cyclic permutations and $E(\sigma)$ is defined by (\ref{E}).
In the second equality $\tau_{1,2}$ denotes the transposition $(1 \ 2)$.
Notice that $E(\sigma)$ appears in both formulas above because, by definition, the
right action of $S_n$ on ${\rm Hom}(V^{\otimes n},V)$ is given by composition with
$E(\sigma)$. Given homogeneous elements $x, y, z \in V$, identities
(\ref{first}) and (\ref{second}) are expressed by:
\begin{gather*}
(-1)^{|x||z|}l_2(l_2(x,y),z) + (-1)^{|x||y|}l_2(l_2(y,z),x)
    + (-1)^{|y||z|}l_2(l_2(z,x),y) = 0  \\[1ex]
\displaystyle{ l_2 (x, y \cdot z) = l_2 (x,y) \cdot z +
                                          (-1)^{(|x|-1)|y|} y \cdot l_2 (x,z) }
\end{gather*}
where the dot product denotes $m_2$.
Notice that the sign $(-1)^{(|x|-1)|y|}$ occurs in the second expression as a consequence
of the transposition $\tau_{1,2}$ and the Koszul sign convention.

To see that an algebra over $H_\bullet (\mathcal{\mathcal{D}})$ is equivalent to
a Gerstenhaber algebra, we just need to define the bracket:
\begin{equation} [x,y] := (-1)^{|x|} l_2(x,y) \label{gbracket}\end{equation}
\noindent
it is not dificult to see that
\begin{equation*}
     [x,y] = - (-1)^{(|x|-1)(|y|-1)}[y,x].
\end{equation*}
and
\[ (-1)^{(|x|-1)(|z|-1)}[[x,y],z] + (-1)^{(|x|-1)(|y|-1)}[[y,z],x] +
                                            (-1)^{(|y|-1)(|z|-1)}[[z,x],y] = 0. \]

This shows that the structure of an algebra over $H_\bullet (\mathcal{\mathcal{D}})$ is equivalent to the structure of
an Gertenhaber algebra on $V$ with bracket defined by (\ref{gbracket}).

\section{The Homology swiss-cheese operad} \label{homology_sc}
In this section we recall the definition of the swiss-cheese operad \cite{Voro99},
denoted by $\mathcal{SC}$. We will present the homology swiss-cheese operad using 2-coloured trees.
Harrelson \cite{Harr04} has presented similarly the homology of open-closed strings in the wider
context of PROPs. The following presentation of the homology swiss-cheese operad is a particular
case of Harrelson's open-closed homology PROP.

The swiss-cheese operad $\mathcal{SC}$ is a 2-coloured operad. We will use the initials of open and closed
as our set of colours: $C = \{ o,c \}$. For $m \geqs 0, n \geqs 0$ with $m + n \geqs 1$,
$\mathcal{SC}(m,n;o)$ is the configuration space of non-overlapping disks labeled $1$ through $m$ and upper
semi-disks labeled $1$ through $n$ embedded by translations and dilations in the
standard unit upper semi-disk so that the embedded semi-disks are all centered on
the diameter of the big semi-disk.

For $m \geqs 1$ and $n = 0$, $\mathcal{SC}(m,0;c)=\mathcal{D}(m)$ is just the usual
component of the little disks operad, and $\mathcal{SC}(m,n;c)$ is the empty set for $n \geqs 1$.
\begin{obs}
The components of the form $\mathcal{SC}(m,0;o)$ were excluded in the original
definition of the operad $\mathcal{SC}$, (see \cite{Voro99}), i.e.,
those components which have only discs as inputs and intervals as output were not
to be considered in  the original definition.
Here, however, we shall keep the components $\mathcal{SC}(m,0;o)$, $m \geqs 1$,
since they are {\it crucial} for the OCHA structure. In fact, as we will see
in this paper, a zero dimensonal generator of the homology of $\mathcal{SC}(1,0;o)$
corresponds to the map $n_{1,0} : L \rar A$.
The physical meaning of $n_{1,0}$ being given by the
``opening of a closed string into an open one'', see \cite{KS06a, KS06b}.
\end{obs}
Let us now describe the homology swiss-cheese operad $H_\bullet (\mathcal{SC})$ in terms of
generators and relations using trees. Since $H_\bullet (\mathcal{SC})$ is a 2-coloured operad,
our trees must also be 2-coloured. The colours we use are wiggly and straight.

Let $l_2 = \!\!\raisebox{-1.5em}{\begin{picture}(0,0)%
\includegraphics{sl2_tree.pstex}%
\end{picture}%
\setlength{\unitlength}{3947sp}%
\begingroup\makeatletter\ifx\SetFigFont\undefined%
\gdef\SetFigFont#1#2#3#4#5{%
  \reset@font\fontsize{#1}{#2pt}%
  \fontfamily{#3}\fontseries{#4}\fontshape{#5}%
  \selectfont}%
\fi\endgroup%
\begin{picture}(665,702)(636,-241)
\put(626,380){\makebox(0,0)[lb]{\smash{{\SetFigFont{8}{9.6}{\rmdefault}{\mddefault}{\updefault}$1$}}}}
\put(1107,380){\makebox(0,0)[lb]{\smash{{\SetFigFont{8}{9.6}{\rmdefault}{\mddefault}{\updefault}$2$}}}}
\end{picture}%
}$ and
$m_2 = \!\!\raisebox{-1.5em}{\begin{picture}(0,0)%
\includegraphics{sm2_tree.pstex}%
\end{picture}%
\setlength{\unitlength}{3947sp}%
\begingroup\makeatletter\ifx\SetFigFont\undefined%
\gdef\SetFigFont#1#2#3#4#5{%
  \reset@font\fontsize{#1}{#2pt}%
  \fontfamily{#3}\fontseries{#4}\fontshape{#5}%
  \selectfont}%
\fi\endgroup%
\begin{picture}(665,702)(636,-241)
\put(626,380){\makebox(0,0)[lb]{\smash{{\SetFigFont{8}{9.6}{\rmdefault}{\mddefault}{\updefault}$1$}}}}
\put(1107,380){\makebox(0,0)[lb]{\smash{{\SetFigFont{8}{9.6}{\rmdefault}{\mddefault}{\updefault}$2$}}}}
\end{picture}%
}$ denote the generators of
$H_\bullet (\mathcal{SC}(2,0;c)) = H_\bullet (\mathcal{D}(2))$.
Both spaces
$\mathcal{SC}(1,0;o)$ and  $\mathcal{SC}(0,2;o)$ are contractible because the elements
of $\mathcal{SC}(1,0;o)$ have only one interior disc while $\mathcal{SC}(0,2;o)$ is homeomorphic
to $\mathcal{C}_1(2)$ (where $\mathcal{C}_1$ is the little intervals operad) which is well known to
be contractible. Their degree zero homology generators will be denoted respectively by
\[ n_{1,0} = \;\raisebox{-1.7em}{\begin{picture}(0,0)%
\includegraphics{n_10_tree.pstex}%
\end{picture}%
\setlength{\unitlength}{3947sp}%
\begingroup\makeatletter\ifx\SetFigFont\undefined%
\gdef\SetFigFont#1#2#3#4#5{%
  \reset@font\fontsize{#1}{#2pt}%
  \fontfamily{#3}\fontseries{#4}\fontshape{#5}%
  \selectfont}%
\fi\endgroup%
\begin{picture}(197,766)(3111,-1305)
\put(3124,-623){\makebox(0,0)[lb]{\smash{{\SetFigFont{8}{9.6}{\rmdefault}{\mddefault}{\updefault}$1$}}}}
\end{picture}%
} \quad \mbox{  and }
   \quad n_{0,2} =\!\! \raisebox{-1.7em}{\begin{picture}(0,0)%
\includegraphics{n_02_tree.pstex}%
\end{picture}%
\setlength{\unitlength}{3947sp}%
\begingroup\makeatletter\ifx\SetFigFont\undefined%
\gdef\SetFigFont#1#2#3#4#5{%
  \reset@font\fontsize{#1}{#2pt}%
  \fontfamily{#3}\fontseries{#4}\fontshape{#5}%
  \selectfont}%
\fi\endgroup%
\begin{picture}(690,743)(3912,-1305)
\put(3912,-646){\makebox(0,0)[lb]{\smash{{\SetFigFont{8}{9.6}{\rmdefault}{\mddefault}{\updefault}$1$}}}}
\put(4418,-655){\makebox(0,0)[lb]{\smash{{\SetFigFont{8}{9.6}{\rmdefault}{\mddefault}{\updefault}$2$}}}}
\end{picture}%
}. \]
So, the degrees of the above generators are:
$|l_2| = 1$ and $|m_2| = |n_{1,0}| = |n_{0,2}| = 0$. The analog of
Cohen's theorem can be stated as follows.
\begin{thm} \label{SC_homology}
The homology swiss-cheese operad $H_\bullet (\mathcal{SC})$ is isomorphic, as a $\mathbb{Z}$-graded
$2$-coloured operad, to the $2$-coloured operad generated by:
$l_2$, $m_2$, $n_{1,0}$ and $n_{0,2}$, satisfying the following relations:
\begin{enumerate}[a)]
  \item $l_2$ is invariant under permutation and satisfies the Jacobi identity;
  \item $m_2$ is invariant under permutation and associative;
  \item $l_2$ and $m_2$ satisfy the Leibniz rule;
  \item $n_{0,2}$ is associative;
  \item %$n_{1,0}$ goes into the center of $n_{0,2}$:
                   \raisebox{-2em}{\input{n_1,0n_0,2.pstex_t}} \hspace{1em} and \hspace{1em}
    %$m_2$, $n_{1,0}$ and $n_{0,2}$ are related by:
              \raisebox{-2em}{\input{1002_tree.pstex_t}}.
\end{enumerate}
\end{thm}
\begin{obs}
  Given a coloured tree, if it has $k$ leaves of some colour $c$, then those
  leaves are labeled $1$ to $k$. So, on the same tree we may have two leaves of
  different colours with the same label.
\end{obs}
\begin{proof}[Proof of Theorem \ref{SC_homology}]
   We first show that $l_2$, $m_2$, $n_{1,0}$ and $n_{0,2}$ generate the operad
   $H_\bullet (\mathcal{SC})$. In fact, $\mathcal{SC}(0,n;o)$ is homeomorphic to the little
   intervals operad $\mathcal{C}_1$ and hence $H_\bullet (\mathcal{SC}(0,n;o))$ is generated
   by $n_{0,2}$. The space $\mathcal{SC}(1,0;o)$ consists of configurations of one disk inside
   the standard semi disk and is also contractible, hence the need for the zero
   dimensional generator $n_{1,0}$. Finally observe that $\mathcal{SC}(m,n;o)$
   is homotopy equivalent to $\mathcal{D}(m) \times S_n$, thus from Theorem \ref{cohen76} any
   class in $H_\bullet (\mathcal{SC}(m,n;o))$ is obtained from operadic composition of
   $l_2$, $n_{1,0}$ and $n_{0,2}$. In order to generate the full operad
   $H_\bullet (\mathcal{SC})$ we need to consider $\mathcal{SC}(m,0;c) = \mathcal{D}(m)$. But we already
   know that its homology is generated by $l_2$ and $m_2$ from Theorem \ref{cohen76}.

   We now show that the generators in fact satisfy the above relations.
   Relation {\it e)} involve only zero dimensional homology classes and
   one can easily check that the compositions indicated in {\it e)} belong to the
   same path component of the swiss-cheese.
   Item {\it d)} follows immediately from the fact that $\mathcal{SC}(0,n;o)$
   is the little intervals operad $\mathcal{C}_1(n)$. Since $\mathcal{SC}(m,0;c) =
   \mathcal{D}(m)$, relations {\it a), b)} and {\it c)} are precisely the statement
   of Theorem \ref{cohen76}.
\end{proof}
We will now study algebras over $H_\bullet (\mathcal{SC})$. Since our main interest is
in OCHA and, as said in the introduction, OCHA is related to {\it part} of the
structure of the homology swiss-cheese, let us define a suboperad
of $H_\bullet (\mathcal{SC})$ containing the relevant structure.
\begin{defi} \label{oc-algebra}
 $\mathcal{OC}$ is the suboperad of $H_\bullet (\mathcal{SC})$ generated by $l_2$,
 $n_{1,0}$ and $n_{0,2}$. Algebras over $\mathcal{OC}$ will be called
 \emph{open-closed algebras} or simply \emph{$OC$-algebras.}
\end{defi}
\begin{obs}\label{lie}
Let $\mathcal{L}$ be the operad defined by $\mathcal{L}(n) = H_{n-1}(\mathcal{D}(n))$ for $n \geqs 1$, i.e.,
$\mathcal{L}$ is the operad of top dimensional homology classes of the little discs operad.
From Theorem \ref{cohen76}, we see that the operad $\mathcal{L}$ is generated by an element
$l_2 \in H_1(\mathcal{D}(2))$ which is invariant
under the symmetric group action and satisfies the Jacobi identity. Consequently, $l_2$ corresponds to a
degree one graded commutative bilinear operation satisfying the Jacobi identity.
Under operadic desuspension, the new generator will have degree zero,
will be graded anti-commutative and will satisfy the Jacobi identity.
In other words, the operadic desuspension $\mathfrak{s}^{-1}$ transforms
$\mathcal{L}$ into the Lie operad:
$\mathcal{L}ie = \mathfrak{s}^{-1} \mathcal{L}$ (see \cite{MSS02} for the definition of operadic (de)suspension).
In this paper we refer to $\mathcal{L}$ as the \emph{Lie operad} by ``abus de langage''.
\end{obs}

There is an analogy between the operads $\mathcal{OC}$ and $\mathcal{L}$ which is sumarized bellow:
%The generator of $\mathcal{L}ie$ belongs to $H_\bullet (\mathcal{D}(2))$ and
%is precisely the top-dimensional class. The generators of $\mathcal{OC}$ belong to
%$H_\bullet (\mathcal{SC}(2,0;c)$, $H_\bullet (\mathcal{SC}(1,0;o)$, $H_\bullet (\mathcal{SC}(0,2;o)$ and
%are precisely their top-dimensional classes.
%The analogy is summarized bellow:
\begin{gather*}
    \mathcal{OC} \quad \Longleftrightarrow \quad
      \mbox{top-dimensional generators of } H_\bullet (\mathcal{SC}) \\
    \mathcal{L} \quad \Longleftrightarrow \quad
      \mbox{top-dimensional generators of } H_\bullet (\mathcal{D}).
\end{gather*}

An algebra over $\mathcal{OC}$ (or $OC$-algebra) consists of a pair of $\mathbb{Z}$-graded
vector spaces $L$ and $A$ such that $L$ is endowed with a degree one symmetric
operation $l_2 : L \otimes L \rar L$ satisfying the Jacobi identity, $A$ has
a degree zero operation $m_2 : A \otimes A \rar A$ defining a structure of
associative algebra and there is a degree zero linear map
$n_{1,0} : L \rar A$. From the first identity in item e) of Theorem \ref{SC_homology},
it follows that $n_{1,0}$ takes $L$ into the center of $A$.

\section{Open-closed homotopy algebras}\label{OCHA}
OCHA's were originally defined in a particular grading and signs
convention where all multilinear maps have degree one and, after being lifted as a
coderivation $D \in {\rm Coder}(S^c L \otimes T^c A)$, the OCHA axioms are
translated into the single condition: $D^2 = 0$.

In order to study the relation between OCHA and the swiss-cheese operad, we
need a definition where grading and signs are given by the corresponding compactified configuration space.
More specifically, a definition where the degrees are equal to minus the
dimension of the configuration space and the signs in the axioms are chosen so as to make
them compatible with the boundary operator in the first row of the $E^1$ term of the spectral
sequence of the compactified configuration space.
In this section we present the definition in this geometrical setting. It is proven in the
Appendix that both definitions are equivalent.

Let us first recall the definition of SH Lie \cite{LS93} algebras in a grading and signs convention
compatible with its compactified configuration space description (see \cite{Voro93,KSV95}).
\begin{defi}[Strong Homotopy Lie algebra]
A strong homotopy Lie algebra (or $L_\infty$-algebra) is a
$\mathbb{Z}$-graded vector space $V$ endowed
with a collection of graded symmetric $n$-ary brackets $l_n : V^{\otimes n} \rar V$,
of degree $3 - 2n$ such that $l_1^{^2} = 0$ and for $n \geqs 2$:
\begin{equation}\label{homotopy_lie}
  \partial \,l_n (v_1, \dots, v_n)\; = \!\!\!\!\! \sum_{\begin{array}{c}\scriptstyle
                                                        \sigma \in \Sigma_{k+l=n} \\[-1ex]
                                                        \scriptscriptstyle{k \geqs 2, l \geqs 1}
                                                      \end{array}}
  \!\!\!\! \epsilon(\sigma)\;
l_{1+l}(l_k(v_{\sigma(1)}, \dots,v_{\sigma(k)}),v_{\sigma(k+1)}, \dots,v_{\sigma(n)}) = 0
\end{equation}
where $\sigma$ runs over all $(k,l)$-unshuffles, i.e., permutations $\sigma \in S_n$
such that $\sigma(i) < \sigma(j)$ for $1 \leqs i < j \leqs k$ and for
$k+1 \leqs i < j \leqs k+l$.
\end{defi}
\begin{obs*}
The operator $\partial$ in the above definition denotes the induced differential on the
endomorphism complex, i.e.:
\[ \partial \,l_n = l_1 l_n + l_n
(l_1 \otimes {\bf 1} \otimes \cdots \otimes {\bf 1} + \; \cdots \; + {\bf 1} \otimes \cdots \otimes
             {\bf 1} \otimes l_1).  \]
\end{obs*}
\begin{defi}[Open-Closed Homotopy Algebra $-$ OCHA] \label{OCHA_def}
An OCHA consists of a 4-tuple $(L,A,\mathfrak{l},\mathfrak{n})$ where $L$ and $A$ are
$\mathbb{Z}$-graded vector spaces,
$\mathfrak{l} = \{ l_n : L^{\otimes n} \rar L \}_{n \geqs 1}$ and
$\mathfrak{n} = \{ n_{p,q} : L^{\otimes p} \otimes A^{\otimes q} \rar A \}_{p+q \geqs 1}$
are two families of multilinear maps where $l_n$ has degree $3 - 2n$ and
$n_{p,q}$ has degree $2 - 2p - q$, such that $(L,\mathfrak{l})$ is an $L_\infty$-algebra
and the two families satisfy the following compatibility condition:
\begin{multline*}
\partial \;n_{n,m}(v_1, \dots, v_n, a_1, \dots, a_m) = \\
= \mquad \sum_{\scriptscriptstyle{\sigma \in \Sigma_{p+r=n}, \;p \geqs 2}}
 \mquad (-1)^{\epsilon(\sigma)}
n_{1+r,m}(l_p(v_{\sigma(1)}, \dots , v_{\sigma(p)}),v_{\sigma(p+1)}, \dots ,
                                          v_{\sigma(n)},a_1, \dots, a_m)\; + \\
%\hspace*{-1.3em}
 + \mquad \mquad  \mquad \!\!\; \sum_{\ \ \ \ \begin{array}{c}
                            \scriptscriptstyle{\sigma \in \Sigma_{p+r=n},\; i+j=m-s}
                            \\[-1ex]
                            \scriptscriptstyle{(r,s) \neq (0,1), (n,m)}
                          \end{array}} \mquad \mquad \mquad \!\!\!\!
                        (-1)^{\mu_{p,i}(\sigma)}
                        n_{p,i+1+j}(v_{\sigma(1)}, . . , v_{\sigma(p)},a_1,
            . . , a_i,
 n_{r,s}(v_{\sigma(p+1)}, . . , v_{\sigma(n)},a_{i+1}, . . , a_{i+s}), a_{i+s+1}, . . , a_m).
\end{multline*}
\begin{multline*}
\mbox{ where } \mu_{p,i}(\sigma) = s + i + si + ms +\epsilon(\sigma) +
s(v_{\sigma(1)} + \cdots + v_{\sigma(p)} + a_1 + \cdots + a_i) + \\ +
(a_1 + \cdots + a_i)(v_{\sigma(i+1)}) + \cdots + v_{\sigma(n)}).
\end{multline*}
\end{defi}
\begin{obs*}
The operator $\partial$ in the above definition denotes the induced differential on the
endomorphism complex, i.e.:
\[ \partial \; n_{n,m} = n_{0,1} n_{n,m} \; - (-1)^m\;
n_{n,m} (d_{L^n} \otimes {\bf 1}_A^{\otimes m} +
{\bf 1}_L^{\otimes n} \otimes d_{A^m}) \]
where
$d_{L^n} =\; \scriptstyle{l_1 \otimes {\bf 1} \otimes \cdots \otimes {\bf 1} + \; \cdots \;
+ {\bf 1} \otimes \cdots \otimes {\bf 1} \otimes l_1}$ and
$d_{A^m} =\; \scriptstyle{n_{0,1} \otimes {\bf 1} \otimes \cdots \otimes {\bf 1} + \; \cdots \; +
{\bf 1} \otimes \cdots \otimes {\bf 1} \otimes n_{0,1}}$.
\end{obs*}

It is convenient to have a shorthand expression for the OCHA relations:
\begin{multline}\label{geo_ocha_shorthand}
\partial \; n_{n,m}  = \mquad \sum_{\scriptscriptstyle{\sigma \in \Sigma_{p+r=n}, \;p \geqs 2}}
n_{1+r,m}(l_p \otimes {\bf 1}_L^{\otimes \,r} \otimes {\bf 1}_A^{\otimes \,m})
(E(\sigma) \otimes {\bf 1}_A^{\otimes \,m})\; +  \\
\hspace*{-1.3em}
 + \mquad \mquad \; \sum_{\begin{array}{c}
                            \scriptscriptstyle{\sigma \in \Sigma_{p+r=n},\; i+j=m-s}
                            \\[-1ex]
                            \scriptscriptstyle{(r,s) \neq (0,1), (n,m)}
                          \end{array}} \mquad \mquad
                        (-1)^{s + i + si + ms}
                        n_{p,i+1+j}({\bf 1}_L^{\otimes \,p} \otimes {\bf 1}_A^{\otimes \,i}
  \otimes n_{r,s} \otimes {\bf 1}_A^{\otimes \,j})(E(\sigma) \otimes {\bf 1}_A^{\otimes \,m})
\end{multline}
\noindent
where $E(\sigma)$ was defined by formula (\ref{E}) on page \pageref{E}.
The complicated sign of the definition is absorbed in the above expression if we assume the
following standard convention:
given two maps $h_1,h_2 : V \otimes W \rar U$, the tensor
product $h_1 \otimes h_2$ defined on $V^{\otimes^2} \otimes W^{\otimes^2}$ is given by:
$(h_1 \otimes h_2)((v_1 \otimes v_2) \otimes (w_1 \otimes w_2)) =
(-1)^{|v_2||w_1|}h_1(v_1,w_1) \otimes h_2(v_2,w_2)$.
\begin{exm}
Here is a list of the first few OCHA relations:
\begin{flalign}
 \ \partial \, n_{0,1} & = 2\,(n_{0,1})^2 = 0 & &
 \label{1st} \\
 \ \partial \, n_{1,1} & = n_{0,2}(n_{1,0} \otimes \id_A) - n_{0,2}(\id_A \otimes n_{1,0}) & &
 \label{2nd} \\
 \ \partial \, n_{2,0} & =
 n_{1,0} l_2 + n_{1,1}(\id_L \otimes n_{1,0}) +
                                     n_{1,1}(\id_L \otimes n_{1,0})E(\tau_{1,2}) & &
 \label{3rd}
\end{flalign}
\begin{multline}
 \label{4th}
 \partial \, n_{1,2} = n_{1,1}({\bf 1}_L \otimes n_{0,2}) - n_{0,2}(n_{1,1} \otimes {\bf 1}_A)
   - n_{0,2}({\bf 1}_A \otimes n_{1,1}) + \\
   + n_{0,3}(n_{1,0} \otimes {\bf 1}_A \otimes {\bf 1}_A)
   - n_{0,3}({\bf 1}_A \otimes n_{1,0} \otimes {\bf 1}_A)
   + n_{0,3}({\bf 1}_A \otimes {\bf 1}_A \otimes n_{1,0})
\end{multline}
\begin{multline}
\label{5th}
 \partial \, n_{2,1} = n_{1,1}(l_2 \otimes {\bf 1}_A) + n_{1,1}({\bf 1}_L \otimes n_{1,1})
   + n_{1,1}({\bf 1}_L \otimes n_{1,1})(E(1\; 2) \otimes {\bf 1}_A) + \\
   + n_{0,2}(n_{2,0} \otimes {\bf 1}_A)
   - n_{0,2}({\bf 1}_A \otimes n_{2,0})
   + n_{1,2}({\bf 1}_L \otimes n_{1,0} \otimes {\bf 1}_A)
   - n_{1,2}({\bf 1}_L \otimes {\bf 1}_A \otimes n_{1,0})\, + \\
   + n_{1,2}({\bf 1}_L \otimes n_{1,0} \otimes {\bf 1}_A)(E(1\; 2) \otimes {\bf 1}_A)
   - n_{1,2}({\bf 1}_L \otimes {\bf 1}_A \otimes n_{1,0})(E(1\; 2) \otimes {\bf 1}_A)
\end{multline}
\end{exm}

Relation (\ref{1st}) simply says that $n_{0,1}$ is a differential operator.
On the other hand, (\ref{2nd}) means that $n_{1,0}$ takes $L$ into the homotopy center of $A$ where
$n_{1,1}$ is the homotopy operator. The configuration space corresponding to $n_{1,1}$ is the cyclohedron $W_2$
(see example \ref{cyclohedra}). The configuration space corresponding to relation (\ref{3rd}) is ``The Eye''
(Figure \ref{theeye_figure} pg.\pageref{theeye_figure}). Relation (\ref{4th}) corresponds to the configuration space $W_3$ (Figure \ref{w3} pg.\pageref{w3}). Finally, relation (\ref{5th}) corresponds to the configuration space illustrated by Figure \ref{c21} pg.\pageref{c21}.
If we consider an OCHA structure where the maps $n_{1,0}$ and
$n_{2,0}$ are set equal to zero, then relations (\ref{4th}) and (\ref{5th}) together say that $n_{1,1} : L \otimes A \rar A$ is a Lie algebra action by derivations up to homotopy.

It is a well known fact that $A_\infty$ and $L_\infty$ algebras can be defined both
in the {\it geometrical setting} (where the degrees of the multilinear maps are minus
the dimension of the corresponding configuration space) and the {\it algebraic setting} where
all the maps have degree 1. It is also well known that both definitions are equivalent
through the (de)suspension operator. The same is true for OCHA.
\begin{prop} \label{coder_ocha}
An OCHA structure $(L,A,\mathfrak{l},\mathfrak{n})$, in the grading and signs conventions of
defintion \ref{OCHA_def}, is equivalent to a degree one coderivation
$D \in {\rm Coder}(S^c (\downarrow \! \downarrow \!\! L) \otimes T^c ( \downarrow \!\! A))$ such that $D^2 = 0$.
\end{prop}

The proof of this fact amounts to an appropriate
use of the Lada-Markl notation for the suspension and desuspension operators $\uparrow$
and $\downarrow$ (see \cite{LM95}) and is provided in the Appendix.

\section{The OCHA operad $\mathcal{OC}_\infty$}\label{OCHA_operad}
\label{ocha_operad}
In this section
we study the operad $\mathcal{OC}_\infty$ whose algebras are precisely OCHA's as given in Definition
\ref{OCHA_def}. Our presentation of $\mathcal{OC}_\infty$ is slightly different from
(but naturally equivalent to) the original definition in \cite{KS06a} because of the different conventions.
We begin by defining the partially planar trees.
\begin{defi} \label{ppt}
  A \emph{partially planar tree} is an isotopy class of oriented rooted trees embedded in the euclidean
  3 dimensional space $\mathbb{R}^3$ such that a fixed subset of edges is contained in the xy-plane.
  Planar edges will be denoted by straight lines, while spatial edges will be denoted by wiggly lines.
\end{defi}
\begin{obs}
Partially planar trees have appeared in the work of Merkulov \cite{Mer02}. Merkulov, however,
uses wiggly lines for planar edges and straight lines for spatial edges.
\end{obs}

The partially planar trees relevant for the definition of $\mathcal{OC}_\infty$ have a specific
form we now begin to describe. We define $l_n$ as the corolla which has
$n$ leaves and only spatial edges and $n_{p,q}$ as the corolla with planar root,
$p$ spatial leaves and $q$ planar leaves. Leaves of different colours are labelled by different sets:
\begin{equation}\label{corollae}
 \raisebox{-3em}{\input{lk.pstex_t} \quad \input{nkl.pstex_t}.}
\end{equation}

As mentioned in the introduction, the grafting operation of a tree $T_2$ on some leaf of a
tree $T_1$ is only defined when the colour of the root of $T_2$ is equal to the colour of
the corresponding leaf of $T_1$. The grafting of a tree $T_2$ with spatial root on
the $i$th spatial leaf of some tree $T_1$ will be denoted by:
\begin{equation}\label{spatial_comp}
 T_1 \circ_i T_2
\end{equation}
On the other hand, the grafting of a tree with planar
root $T_4$ on the $i$th planar leaf of some tree
$T_3$ will be denoted by:
\begin{equation}
 T_3 \bullet_i T_4
\end{equation}
Consider the set of all corollas $n_{p,q}$ and $l_n$ with $2p + q \geqs 2$ and $n \geqs 2$.
Let $\mathcal{T}(n)$ denote the set of all partially planar trees $T$ with
$n$ leaves which can be obtained by grafting a finite number of corollas in the
above set.
Let $\mathcal{T}_o(p,q) \subseteq \mathcal{T}(p+q)$
denote the subset of trees with planar root having $p$
spatial leaves and $q$ planar leaves.
Let $\mathcal{T}_c(n) \subseteq \mathcal{T}(n)$ be the subset
of trees with spatial root.
\begin{defi}
For $2p + q \geqs 2$, we define $\mathcal{N}_{\infty}(p,q)$ as the vector space spanned by
$\mathcal{T}_o(p,q)$ and for $n \geqs 2$, $\mathcal{L}_{\infty}(n)$ is defined as the vector space spanned by $\mathcal{T}_c(n)$. The space $\mathcal{N}_{\infty}(0,1)$ is defined as the vector space spanned by the tree
with only one planar edge and no vertices, while $\mathcal{L}_{\infty}(1)$ is defined similarly as the
vector space spanned by the tree with only one spatial edge and no vertices.
\end{defi}
\begin{obs}
Notice that if a tree in $\mathcal{T}(n)$ has a spatial root, then all of its edges must also be spatial
because of the corollas we have chosen as generators.
\end{obs}
\noindent
Let $|i(T)|$ be the number of internal edges of $T$,
we define the degree of $T \in \mathcal{T}(n)$ as follows:
\begin{equation}\label{degreeofT}
|T| = \left\{ \begin{array}{ll}
                  |i(T)| + 2 -2p - q, & {\rm if}\; T \in \mathcal{T}_o(p,q) \\[1ex]
                  |i(T)| + 3 - 2n,    & {\rm if}\; T \in \mathcal{T}_c(n)
                 \end{array} \right.
\end{equation}
\noindent
in particular, $|n_{p,q}| = 2 - 2p - q$ and $|l_{n}| = 3 - 2n$.
Now we can define the spaces: $\mathcal{N}_{\infty} = \bigoplus_{k+l \geqs 1} \mathcal{N}(k,l)$ and
$\mathcal{L}_{\infty} = \bigoplus_{n \geqs 1} \mathcal{L}_\infty(n)$, and finally define:
\begin{equation}
  \mathcal{OC}_\infty = \mathcal{L}_\infty \oplus \mathcal{N}_\infty.
\end{equation}

There is a symmetric group action on spatial leaves by permuting the labels of the spatial leaves,
and there is no symmetric group action on planar leaves. In other words, given a tree $T \in \mathcal{OC}_\infty$
with $p$ spatial leaves and $q$ planar leaves, the group $S_p$ acts on $T$ by permuting the
labels of the spatial leaves, while the planar leaves remain fixed.

The space $\mathcal{OC}_\infty$ we have just defined has the structure of a 2-coloured operad of
graded vector spaces defined by the grafting operations $\circ_i$ and $\bullet_i$ and by the
symmetric group action on spatial leaves. Let us now define a differential operator
$d: \mathcal{OC}_\infty \rar \mathcal{OC}_\infty$. We proceed analogously to the definition
given in \cite{HS93} (see also \cite{KSV95}).

Let us first define the action of $d$ on corollas $l_n$ and $n_{p,q}$:
\begin{equation}\label{dl}
  d \ l_n = \mquad\sum_{\begin{array}{c}
                    \scriptstyle k+l=n+1
                                 \\[-1ex]
                    \scriptstyle k,l \geqs 2
                  \end{array}}\mquad
            \sum_{\begin{array}{c}
                    \scriptstyle  {\rm unshuffles}\; \sigma:  \\[-.8ex]
                    \scriptstyle \{ 1,2,\dots,n \} = I_1 \cup I_2 \\[-.8ex]
                    \scriptstyle \#I_1 = k, \ \#I_2 = l - 1
                  \end{array}} \raisebox{-4em}{\input{dif_lieII.pstex_t}}
\end{equation}
observing that an unshuffle $\sigma$ is equivalent to a partition $(1,2, \dots, n) = I_1 \cup I_2$ into two ordered subsets $I_1$ and $I_2$. On the other hand:
$d \; n_{n,m} = $
\begin{equation}\label{dn}
   \hskip -1em = \mquad\sum_{\begin{array}{c}
                                  \scriptstyle k+l=n+1
                                  \\[-1ex]
                                  \scriptstyle k,l \geqs 2
                               \end{array}}\mquad
                         \sum_{\begin{array}{c}
                                  \scriptstyle  {\rm unshuffles}\; \sigma:  \\[-.8ex]
                                  \scriptstyle \{ 1,2,\dots,n \} = I_1 \cup I_2 \\[-.8ex]
                                  \scriptstyle \#I_1 = k, \ \#I_2 = l - 1
                               \end{array}} \mquad \Bigg( \ \raisebox{-4em}{\input{dif_ocha1.pstex_t}} +
                   \sum_{0 \leqs i,s \leqs m} (-1)^{s+i+si+ms}\raisebox{-4em}{\input{dif_ocha2.pstex_t}}\; \Bigg).
\end{equation}
Once $d$ is defined on the generators of $\mathcal{OC}_\infty$, it is extended to the whole operad
by the Leibniz rule:
\[ d \, (T \circ_i T_1) = d \, T \circ_i T_1 + (-1)^{|T|} T \circ_i d \, T_1
\qquad d \, (T \bullet_i T_2) = d \, T \bullet_i T_2 + (-1)^{|T|} T \ \bullet_i d \, T_2 \]
where: $T_1$ is a tree with spatial root and $T_2$ is a tree with planar root.
With the operator $d$, $\mathcal{OC}_\infty$ becomes a differential graded $2$-coloured operad.
\begin{obs}\label{t'->t}
For trees in $\mathcal{OC}_\infty$, let $T^{'} \rar T$ indicate that $T$ is obtained from $T^{'}$ by contracting a
spatial or planar internal edge. The above defined differential operator
$d : \mathcal{OC}_\infty \rar \mathcal{OC}_\infty$ is, up to sign, simply given by:
\[ d(T) = \sum_{T^{'} \rar T} \pm\; T^{'}, \]
the only difference between the above definition and the original one in \cite{KS06a} is the sign.

According to the grading defined by (\ref{degreeofT}), the operator $d$ has degree $1$. So $\mathcal{OC}_\infty$ is
an operad of cochain complexes.
\end{obs}
Given two differential graded vector spaces $L$ and $A$, we say that $(L,A)$ is an algebra over
$\mathcal{OC}_\infty$ if there is a morphism of differential graded 2-coloured operads:
\[ \Psi : \mathcal{OC}_\infty \rar {\rm End}_{L,A}    \]
where ${\rm End}_{L,A}$ is the 2-coloured endomorphism operad of the pair $(L,A)$, as described by (\ref{end_C}).
Since $\Psi$ is a chain map and the differential operator $\partial$ on ${\rm End}_{L,A}$
is precisely the one used in formulas (\ref{homotopy_lie}) and (\ref{geo_ocha_shorthand}),
it follows that $(L,A)$ is an algebra over $\mathcal{OC}_\infty$ if, and only if, it admits
the structure of an OCHA.
\begin{obs}\label{ideal_N}
By definition, $\mathcal{OC}_\infty = \mathcal{L}_\infty \oplus \mathcal{N}_\infty$, where
$\mathcal{N}_\infty$ is spanned by trees with planar root. For trees in $\mathcal{OC}_\infty$, grafting
two trees $T_1$ and $T_2$ where at least one of them has a planar root always results in a tree with
a planar root. So, $\mathcal{N}_\infty$ is an ideal in $\mathcal{OC}_\infty$.
On the other hand, $\mathcal{L}_\infty$ is a suboperad of $\mathcal{OC}_\infty$, since trees with spatial
root can only have spatial edges. Finally, we observe that $\mathcal{OC}_\infty$ has a strucutre of
module over $\mathcal{L}_\infty$ given by the grafting operation in $\mathcal{OC}_\infty$.
\end{obs}

\section{The compactification $\overline{C(p,q)}$}\label{c(p,q)}
In this section we recall the construction of the space $\overline{C(p,q)}$, first introduced by
Kontsevich in \cite{Kont03}. We use the fact that $\overline{C(p,q)}$ is a manifold with corners
and study the combinatorics of its boundary strata to show that
the first row of the $E^1$ term of the spectral sequence determined by $\overline{C(p,q)}$
is isomorphic, as a differential complex, to $\mathcal{N}_\infty(p,q)$.

Let $p,q$ be non-negative integers satisfying the
inequality $2p+q \geqs 2$. We denote by Conf$(p,q)$ the configuration space of marked points
on the upper half plane $H = \{ z \in \mathbb{C} \;|\; {\rm Im}(z) \geqs 0 \}$ with
$p$ points in the interior and $q$ points on the boundary (real line):
\begin{multline*}
    {\rm Conf}(p,q) = \{ (z_1, \dots, z_p, x_1, \dots, x_q) \in H^{p+q} \;|\;
   z_{i_1} \neq z_{i_2},\, x_{j_1} \neq x_{j_2} \; \forall i_1 \neq i_2, j_1 \neq j_2 \\
   {\rm Im}(z_i) > 0,\,  {\rm Im}(x_j) = 0 \; \forall i, j \}
\end{multline*}

The above configuration space ${\rm Conf}(p,q)$ is the cartesian product of an open subset of $H^p$ and
an open subset of $\mathbb{R}^q$ and, consequently, is a $2p + q$
dimensional smooth manifold.
Let $C(p,q)$ be the quotient of ${\rm Conf}(p,q)$ by the action of the
group of orientation preserving affine transformations that leaves the real line fixed:
\[ C(p,q) = {\rm Conf}(p,q)\Big/(z \mapsto az + b) \quad
a,b \in \mathbb{R}, \; a > 0 . \]
The condition $2p+q \geqs 2$ ensures that the action
is free and thus $C(p,q)$ is a $2p + q - 2$ dimensional smooth manifold.

Let ${\rm Conf}_n(\mathbb{C})$ be the configuration space
of $n$ points in the complex plane.
We take the quotient by affine transformations
$z \mapsto az + b$ where $a \in \mathbb{R}$, $a > 0$ and $b \in \mathbb{C}$ and
define $C(n) := {\rm Conf}_n(\mathbb{C}) \big/ (z \mapsto az + b)$.
Again $C(n)$ is  a smooth manifold. The real version of the Fulton-MacPherson
compactification $\overline{C(n)}$ is defined in the usual way, see \cite{FM94,AS94,
Sinha04}. Let $\phi$ be the embedding:
\begin{equation}\label{embedd} \phi : C(p,q)  \lra C(2p + q)  \end{equation}
defined by $\phi(z_1 , \dots, z_p, x_1, \dots, x_q) =
(z_1 ,\bar{z}_1, \dots, z_p, \bar{z}_p, x_1, \dots, x_q)$, where
$\bar z$ denotes complex conjugation.
\begin{defi}
The compactification of the configuration space $C(p,q)$ is defined as the
closure in $\overline{C(2p + q)}$ of the image of $\phi$.
It will be denoted by $\overline{C(p,q)}$.
\end{defi}

The compactification $\overline{C(n)}$ of points in the plane can
be intuitivelly described through ``{\it bubbling offs}'' on the sphere
(the one point compactification of the plane). In the case of
$\overline{C(p,q)}$, one can think of the closed disc as the one point
compactification of the upper half plane and think of the embedding $\phi$
as taking the closed disc to the upper hemisphere of the above sphere. Punctures in the
bulk of the disc are reflected through the equator.
Points in $\overline{C(p,q)}$ can be intuitivelly described through ``{\it bubbling offs}''
on the disc. Those bubbling offs are pictured in the next figures.

\begin{figure}[h!]
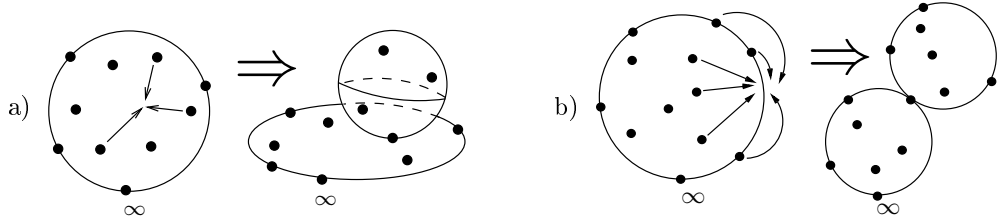

\centering
a) \;\raisebox{-4em}{\input{bubb4.pstex_t}}  \quad\quad\quad
b) \;\raisebox{-4em}{\input{bubb3.pstex_t}}
\caption{The two possible types of bubbling off on the closed disc. }
\label{bulkboundary}
\end{figure}

\subsection{The Stratification of $\overline{C(p,q)}$}\label{stratification}

The combinatorics of the compactification $\overline{C(n)}$ of the configuration space
of points in the complex plane is well known to be described in terms of trees.
In other words, its boundary strata can be labeled by trees
(see: \cite{AS94,FM94,KSV95,Sinha04,Voro00,Voro01}).
Since $\overline{C(p,q)}$ was defined through the embedding $\phi: C(p,q) \rar C(2p+q)$,
it naturally inherits its combinatorics from that of $\overline{C(2p+q)}$. Leaves corresponding
to the $p$ points in the bulk of the upper half plane are spatial, while leaves corresponding
to the $q$ points on the boundary (real line) are planar. The combinatorics of $\overline{C(p,q)}$
is thus described by partially planar trees. We follow the notation of \cite{KSV95} and state this fact
in the following theorem.

\begin{thm}\label{strata}
There is a stratification of $\overline{C(p,q)}$ such that:
\begin{enumerate}[$(1)$]
 \item $ \overline{C(p,q)} = \coprod_{T \in \mathcal{T}_o(p,q)} S_T $.
 Each stratum $S_T$ is a smooth submanifold and
 ${\rm codim}_{\mathbb{R}} S_T = |i(T)| = \mbox{ number of internal edges of $T$ }$;
 \item there is a unique open stratum $S_{n_{p,q}} = C(p,q) \quad 2p + q \geqs 2$;
 \item for each tree $T \in \mathcal{T}_o(p,q)$ we have the identity
       \[ S_T = S_{n_{p_1,q_1}} \times S_{\delta_1} \times \cdots \times
                                                   S_{\delta_n}  \]
where each $\delta_i$ is a corolla of the form $n_{k,l}$ or $l_k$, and $T$ is
obtained by grafting the corollas $\delta_1, \dots, \delta_n$ to $n_{p_1,q_1}$.
 \item The boundary of the closure $\overline{S_T}$ of each stratum is given
by $\partial \overline{S_T} = \bigcup_{T' \rar T} \overline{S_{T'}}$, \\
where $T' \rar T$ means that $T$ is obtained from $T'$ by contracting a
internal edge.
\end{enumerate}
\end{thm}

In case $p=0$, the space $\overline{C(0,q)}$ is the associahedron
$K_q$ \cite{Sta63a,Sta63b}. The labeling
of the boundary strata of $\overline{C(0,q)}$, in this case, reduces to the well
known labeling of the facets of $K_q$ by planar trees.
\begin{exm}[$\overline{C(1,q)}$ is the cyclohedron $W_{q+1}$]
 \label{cyclohedra}
 The cyclohedron  was introduced by Bott and Taubes
 \cite{BT94} and received
 its name from Stasheff \cite{Sta97}. It is defined as the Fulton-MacPherson
 compactification of the configuration space of points on the circle
 $S^1$ modded out by the group of rotations SO$(2) = S^1$. The equivalence between
 $\overline{C(1,q)}$ and $W_{q+1}$ will be described below in terms of partially planar
 trees using the above theorem.

 By fixing the interior point of $C(1,q)$ to be equal to $i \in \mathbb{C}$, the remaining
 points are on the real line. It is not difficult to see that $C(1,q)$ is an open simplex homeomorphic to
 the configuration space of points on the circle modded out by SO$(2)$. The compactification $\overline{C(1,q)}$ is obtained
 from that open simplex by performing iterated blow ups. Hence $\overline{C(1,q)}$ is a polytope.
 In order to show that
 $\overline{C(1,q)}$ and $W_{q+1}$ are equivalent polytopes, we just need to establish
 a one-to-one correspondence between bracketings around the $q+1$ marked points on the circle ($q$ points on the real line
 plus one point marked $\infty$)
 and the partially planar trees in $\mathcal{T}_o(1,q),$ showing also that the correspondence respects
 the incidence relations.

 In fact, for any $q \geqs 0$ the facets of the cyclohedron $W_{q+1}$ are
 labeled by (i.e. are in one-to-one correspondence with)  all the meaningful
 ways of inserting  brackets in an expression of $q + 1$ letters disposed
 on a circle. The codimension of the facet corresponding to a given bracketing
 is equal to the number of brackets inserted, as illustrated by \emph{Figure \ref{w3}}
 (see also Devadoss' paper \cite{Deva01} on the cyclohedra.).

First recall that $\overline{C(1,q)}$ can be described as the compactified
configuration space of points on the closed disc (the one point compactification of the
upper half plane) with $1$ point in the interior of the disc, $q$ points on
the boundary of the disc plus one boundary point marked as $\infty$.
From the ``{\it bubbling off}'' description of points in the compactification, we know that
each facet of codimension $k$ in $\overline{C(1,q)}$ corresponds to $k+1$ discs
joined at points in the boundary such that exactly one of those discs contains one
point in the bulk while the others contain only points in the boundary.

Let us exhibit the correspondence between `circular bracketings' and `bubbling offs'.
Consider a point in $\overline{C(1,q)}$ in the bubbling off description, as a number
of discs joined at ``double points''. There is only one of these discs which contains
the interior point, the remaining ones only contain points in the boundary.
The correspondence goes as follows (see Figure \ref{brabubb}):
\begin{enumerate}
\item the disc containing the interior point corresponds to the circle;
\item points on the boundary of the disc containing the interior point correspond to
points on the circle which are not inside any bracket;
\item a disc joined to the disc containing the interior point corresponds to a
bracketing; two discs joined correspond to a bracketing inside another bracketing
or to two disjoint bracketings, and so on.
\end{enumerate}
In order to get a tree from the bubbling off, we associate to the discs their {\it dual graphs}.
According to the usual procedure,
each disc correspond to a vertex; the point marked $\infty$ corresponds to the root;
the double points correspond to the edges and the remaining marked points correspond
to the leaves. Since the correspondence between bracketings on $S^1$ and joining discs
is established, there follows the correspondence between bracketings on $S^1$ and trees in
$\mathcal{T}_o(1,q)$ (see Figure \ref{brabubb}).
Since the facets of both polytopes are in a one-to-one correspondence compatible with their corresponding
boundaries, it follows that $W_q = \overline{C(1,q)}$.
\begin{figure}[htb]
 \centering
 \hspace*{\stretch{2}}
 \input{bra_bubb.pstex_t}
 \hspace*{\stretch{1}}
 \caption{Example of the correspondence between circular bracketings and trees.}
 \label{brabubb}
\end{figure}
\end{exm}
Figures \ref{w3} and \ref{theeye_figure} on page \pageref{w3}
illustrate the spaces $\overline{C(1,2)}$ and $\overline{C(2,0)}$. Figure \ref{c21} on page \pageref{c21}
is a portrait of $\overline{C(2,1)}$ due to S. Devadoss where
we have included the partially planar trees corresponding to its codimension one boundary strata (see also \cite{KS06c}).
The OCHA relations corresponding to the spaces $\overline{C(2,0)}$, $\overline{C(1,2)}$ and $\overline{C(2,1)}$
are given in formulas (\ref{3rd}), (\ref{4th}) and (\ref{5th}).
\subsection{ The space $\overline{C(p)}$ as a deformation retract of $\overline{C(p,q)}$ } \label{contract}
We close the present section by pointing out a fact that will play a crucial role in the proof of our main theorem
(Theorem \ref{main}).
There is a stratum $S_{T}$ in $\overline{C(p,q)}$ which is homeomorphic to $C(p)$, where $T$ is the following tree:
\begin{equation}\label{T}
 T \quad = \raisebox{-3em}{\input{Cp_copy.pstex_t}},
\end{equation}
moreover, $S_{T}$ is a deformation retract of $\overline{C(p,q)}$.

%$C(p,q) = {\rm Conf}(p,q)\Big/(z \mapsto az + b)$
In fact, by putting a collar neighborhood along the boundary, we see that there is a deformation retraction of $\overline{C(p,q)}$ onto
$C(p,q) = \overline{C(p,q)} \; \backslash \; \partial \,\overline{C(p,q)}$.
This last space was defined as
the quotient of ${\rm Conf}(p,q)$ by the group of affine transformations $z \mapsto az + b$, where
$a,b \in \mathbb{R}$. Since that group is contractible, $C(p,q)$ is
homotopy equivalent to ${\rm Conf}(p,q)$. Now,
${\rm Conf}(p,q)$ is homeomorphic to ${\rm Conf}_{\mathbb{C}}(p) \times {\rm Conf}_{\mathbb{R}}(q)$
and ${\rm Conf}_{\mathbb{R}}(q)$
is well known to be contractible. Thus, by composing all those contractions and homotopy equivalences, we get the
claimed deformation retraction of $\overline{C(p,q)}$ onto $S_T$. To see that $S_T$ is in fact homeomorphic to
$C(p)$, just notice that $T$ is obtained by grafting the trees $l_p$, $n_{1,0}$ and a binary planar tree. The
configuration space corresponding to $n_{1,0}$ and to any binary planar tree is just one point, the homeomorphism thus
follows from Theorem \ref{strata}.
Notice that the retraction is essentialy determined by the contraction of ${\rm Conf}_{\mathbb{R}}(q)$ to a single point.
This means that the configuration of the interior points are unaffected during the contraction of $\overline{C(p,q)}$
onto $\overline{S_T} = \overline{C(p)}$.

Those facts will be used in the next section
along with the fact, proven by P. May in \cite{May72}, that $C(p)$ is $S_p$-equivariantly homotopy equivalent to $\mathcal{D}(p)$, where $\mathcal{D}$ denotes the little discs operad.

\section{OCHA and the spectral sequence of $\overline{C(p,q)}$}\label{OCHA_SS}

In this section we show that the first row of the $E^1$ term of the spectral
sequence of $\overline{C(p,q)}$ is isomorphic, as a chain complex, to $\mathcal{N}_{\infty}(p,q)$. The isomorphism is natural with respect to the operad
composition. This fact depends crucially on the study of the stratification
of $\overline{C(p,q)}$ as a manifold with corners.

%\subsection{The spectral sequence}
Every compact manifold with corners induces a spectral sequence converging to its homology. In fact,
the boundary strata of the manifold induces a natural filtration on its singular chain complex which
ensures the existence of the spectral sequence. Since the boundary filtration is finite, the spectral
sequence is convergent.

Let us study the spectral sequence in the case of the manifold $\overline{C(p,q)}$.
Consider the topological filtration of $\overline{C(p,q)}$:
\begin{equation*}
 F^i \overline{C(p,q)}  = \{ \mbox{closure of the union $\coprod_{T} S_T$
 of strata of dimension $i$} \}
 = \bigcup \{ \overline{S_T} \;|\; {\rm dim}S_T = i \}
\end{equation*}
We will denote $F^i \overline{C(p,q)}$ more simply by $F^i$, with $2p + q - 2 \geqs i \geqs 0$,
remembering that the dimension of $C(p,q)$ is $2p + q - 2$.
The topological filtration induces a filtration on the singular chain complex
of $\overline{C(p,q)}$ and we have the spectral sequence.
\begin{thm}\label{spec_seq}
 There is a spectral sequence $E^r_{m,n}$ converging to $H_*(C(p,q))$. Its $E^1$
term has the form $E^1_{m,n} = H_{m+n}(F^m,F^{m-1})$ and, for $n=0$, the complex
% \begin{enumerate}[(1)]
\[ 0 \rar E^1_{2p+q-2,0} \rar \cdots \rar E^1_{m,0} \rar E^1_{m-1,0} \rar \cdots
                         \rar E^1_{0,0} \rar 0 \]
is isomorphic to the $p,q$ component $\mathcal{N}_\infty(p,q)$ of the ideal
$\mathcal{N}_\infty \lhd \mathcal{OC}_\infty$.
% \end{enumerate}
\end{thm}
%\noindent See Observation \ref{ideal_N} at pg.\pageref{ideal_N}
%for the definition of $\mathcal{N}_\infty$. 
%The above theorem is essentially a consequence
%of Theorem \ref{strata}. In the following proof we just employ a well known argument.
\begin{proof}
From Theorem \ref{strata}, we have:
$F^m \,\backslash\, F^{m-1} = \displaystyle{\coprod_{|T| = - m} \mquad S_T}$,
where $|T|$ is the degree of $T$ defined by (\ref{degreeofT}). The equality $|T| = - m$ means
that the number of internal edges of $T$ is equal to the codimension of its corresponding submanifold
$S_T$. Using the Lefshetz duality theorem:
\begin{equation}\label{first_term}
 E^1_{m,0} = H_m(F^m,F^{m-1})
	   = H^0(F^m \,\backslash\, F^{m-1})
           = H^0(\!\! \coprod_{|T| = - m} \mquad S_T)
           = \bigoplus_{|T| = - m} k.
\end{equation}
As a vector space, $E^1_{m,0}$ is thus exactly the vector space
generated by trees of degree $-m$ in $\mathcal{N}_\infty(p,q)$.
The spectral sequence is a homology spectral sequence, hence $|d^1| = -1$.
On the other hand, $\mathcal{N}_\infty$ is an ideal in $\mathcal{OC}_\infty$
which is an operad of cochain complex. So the diferential $d$ in the ideal
$\mathcal{N}_\infty$ has degree $1$. This is consistent with the fact, used in (\ref{first_term}),
that the degree of a tree is minus the dimension of its corresponding relative homology class.

Now we need to check that the differential $d^1$ on the first row of the
spectral sequence coincides with $d$ defined by formulas (\ref{dl}) and (\ref{dn}).

In fact, by Theorem \ref{strata} and the Lefshetz duality, each relative class in
$H_m(F^m,F^{m-1})$ is given by the closure $\overline{S_T}$ of a submanifold $S_T$.
To see that $d^1$ coincides with $d$, recall that the operator $d^1$ is
given by the relativization of the boundary operator $\partial$:
\begin{diagram}[h=1em,w=2em]
                    &                       & H_{m-1}(F^{m-1}) &                 &                            \\
                    &\ruTTo^{\partial_m}    &                  &\rdTTo^{j_{m-1}} &                            \\
 H_{m}(F^m,F^{m-1}) &                       &\rTTo^{d^1_m}     &                 & H_{m-1}(F^{m-1},F^{m-2}).  \\
\end{diagram}
Item 4 of Thm. \ref{strata}, says that:
$\partial \overline{S_T} = \bigcup_{T' \rar T} \overline{S_{T'}}$.
Observation \ref{t'->t} (pg. \pageref{t'->t}) implies that $d^1 = d$.
\end{proof}
\begin{obs}
It is still necessary to check that the signs given by the coboundary operator in the spectral sequence
coincide with the signs given in formulas (\ref{dl}) and (\ref{dn}). That is a somewhat tedious exercise which
consists of comparing the orientation induced on the product of two oriented manifolds with the orientation induced
by the operadic embedding of that product manifold into the boundary strata of other oriented manifold.
\end{obs}

\subsection{The Quasi-isomorphism of $\mathcal{L}_\infty$-modules} \label{quism}
According to Definition \ref{oc-algebra}, the open-closed operad $\mathcal{OC}$ is the
operad generated by top-dimensional homology classes
of the swiss-cheese operad, i.e., $\mathcal{OC}$ is the suboperad of $H_\bullet (\mathcal{SC})$ generated by
$n_{1,0}$, $l_2$ and $n_{0,2}$ (see Definition \ref{oc-algebra}).
Recall $\mathcal{L}$ is the suboperad of $H_\bullet(\mathcal{D})$ defined by $\mathcal{L}(n) = H_{n-1}(\mathcal{D}(n))$ for $n \geqs 1$.
We refer to $\mathcal{L}$ as the Lie operad, since algebras over it are equivalent to Lie algebras
(see Observation \ref{lie} pg. \pageref{lie}).
From the tree description of $H_\bullet (\mathcal{SC})$ given in section \ref{homology_sc}, we see that
$\mathcal{OC}$ is a suboperad of $\mathcal{OC}_\infty$.

We know that $\mathcal{OC}_\infty$ is an $\mathcal{L}_\infty$-module (see Observation \ref{ideal_N}) and
that $\mathcal{OC}$ is an $\mathcal{L}$-module, since $\mathcal{L}$ is an suboperad of $\mathcal{OC}$.
Consequently, $\mathcal{OC}$ has a natural structure of $\mathcal{L}_\infty$-module induced by
the well known quasi-isomorphism of operads:
$\mu : \mathcal{L}_\infty \rar \mathcal{L}$ defined by $\mu(l_2) = l_2$ and $\mu(l_n) = 0$
for $n \geqs 3$ (see \cite{MSS02,Markl00,Markl99} for details).

Proposition \ref{L_morphism} in the Appendix says that there is a morphism
of differential graded $\mathcal{L}_\infty$-modules extending the identity on $\mathcal{OC}$:
\begin{equation}\label{quism}
\mu : \mathcal{OC}_\infty \rar \mathcal{OC}.
\end{equation}
The $\mathcal{L}_\infty$-morphism $\mu$ vanishes on the corollae that are not in $\mathcal{OC}$.
The restriction of $\mu$ to $\mathcal{L}_\infty$ coincides with the above mentioned quasi-isomorphism
between $\mathcal{L}_\infty$ and $\mathcal{L}$.
Since $\mathcal{OC}_\infty = \mathcal{L}_\infty \oplus \mathcal{N}_\infty$,
in order to prove that $\mu$ is a quasi-isomorphism of $\mathcal{L}_\infty$-modules,
we need to study the cohomology of the ideal $\mathcal{N}_\infty$(see corollary \ref{N_coro}).

Let us begin by showing that, for any $p,q$ such that $2p + q \geqs 2$,
the cohomology of $\mathcal{N}_\infty(p,q)$ is
isomorphic, as $\mathbb{Z}$-graded vector spaces, to $H_\bullet (\mathcal D(p))$
(where $\mathcal D$ is the little discs operad).
From Theorem \ref{spec_seq}, $\mathcal{N}_\infty(p,q)$ is isomorphic to the complex given by:
\[ 0 \lra H_{2p+q-2}(F^{2p+q-2},F^{2p+q-3}) \lra \cdots \lra H_{m}(F^m,F^{m-1})
 \lra \cdots \lra H_0(F^0) \lra 0 \]
where each $F^i$ is the closure of the disjoint union of the $i$-dimensional strata
of $\overline{C(p,q)}$.

Notice that $E^r_{m,n}$ is a homology spectral sequence converging to $H_*(C(p,q))$.
In the proof of Theorem \ref{spec_seq}, we have seen that the above complex is isomorphic to a complex generated by trees whose
degree is minus the degree of their corresponding relative homology classes. So, it is a cochain complex
generated by trees, precisely: $\mathcal{N}_\infty(p,q)$. In what follows, we use this to establish a relation
between the {\it cohomology} of $\mathcal{N}_\infty$ and the {\it homology} of the little disks operad.

Since each stratum is a smooth submanifold, it follows that
$F^i$ has the homotopy type of a CW complex. The manifold $\overline{C(p,q)}$ has thus
the homotopy type of a CW complex $X$ such that each skeleton $X^i$ is homotopy equivalent
to $F^i$. It is well known that for any CW complex, the map $H_n(X^n) \rar H_n(X^{n+1}) \simeq H_n(X)$,
induced by the inclusion $X^n \hookrightarrow X^{n+1}$, is
surjective for all $n$. Consequently, the map
\[ H_n(F^n) \rar H_n(F^{n+1}) \simeq H_n(\overline{C(p,q)}) \]
is also surjective. For the same reason, we have: $H_n(F^{n-1}) \simeq H_n(X^{n-1}) = 0$. Now, consider
the usual commutative diagram:
\begin{diagram}[h=1em,w=2em]
%      &       &            &                         &  &         & &  &  0        &                          & \\
%      &       &            &                         &  &         & &\ruTTo &        &                          & \\
       &       &  0         &                         &  &         & H_n(F^{n+1}) &
     &\mquad \simeq H_n(\overline{C(p,q)}) &          &     &  \\
       &       &            & \rdTTo                  &  & \ruTTo  &  &              &                          & \\
       &       &            &                         & H_n(F^n) & &  &              &                          & \\
       &       &            &\ruTTo^{\partial_{n+1}}  & &\rdTTo^{j_n}& &            &                          & \\
\cdots & \rTTo & H_{n+1}(F^{n+1},F^n)& &\rTTo^{d_{n+1}}& & H_{n}(F^n,F^{n-1})& & \rTTo^{d_{n}}&
& H_{n-1}(F^{n-1},F^{n-2}) & \rTTo  & \cdots,   \\
       &       &            &        & &           & & \rdTTo_{\partial_n} &                & \ruTTo_{j_{n-1}} \\
       &       &            &        & &           & &                     &H_{n-1}(F^{n-1})& \\
       &       &            &        & &           & &  \ruTTo             &                & \\
       &       &            &        & &           &0&                     &                & \\
       &       &            &        & &           & &                     &                & \\
\end{diagram}
since $H_n(F^n) \rar H_n(F^{n+1}) \simeq H_n(\overline{C(p,q)})$ is surjective and
$H_n(F^n) \stackrel{j_n}{\lra} H_n(F^{n},F^{n-1})$ is injective, from the exactness of the sequence of the pair
$(F^n,F^{n-1})$ one can see that the $n$th cohomology group of the complex $\mathcal{N}_{\infty}(p,q)$ is isomorphic to $H_n(\overline{C(p,q)})$. As observed before, $\overline{C(p,q)}$ is homotopy equivalent to $\mathcal{D}(p)$, so the following lemma is proved.
\begin{lem}\label{1st_isomorphism}
 $H^{k}(\mathcal{N}_\infty(p,q)) \simeq H_k(\mathcal{D}(p))$ for every $k \geqs 0$ and $p,q$ such that $2p+q\geqs2$.
\end{lem}
For any $q \geqs 0$, consider the following sequence of vector spaces:
\[ H^\bullet (\mathcal{N}_\infty(\_,q)) := \{ H^\bullet (\mathcal{N}_\infty(p,q)) \}_{p \geqs 1}. \]
Since $\mathcal{L}$ is just the operad generated by a binary tree $l_2$ of degree 1 which is invariant
under the action of the symmetric group $S_2$ and satisfies the Jacobi identity, there is a natural
injection $\mathcal{L} \hookrightarrow H^\bullet (\mathcal{OC}_\infty)$. Since $\mathcal{N}_\infty$ is an ideal
in $\mathcal{OC}_\infty$, it follows that $H^\bullet (\mathcal{N}_\infty)$ is an ideal in $H^\bullet (\mathcal{OC}_\infty)$.
Consequently, for any $q \geqs 0$ we have a structure of $\mathcal{L}$-module on $H^\bullet (\mathcal{N}_\infty(\_,q))$.
Since $\mathcal{L} = \{ H_{n-1}(\mathcal{D}(n)) \}_{n \geqs 1}$ is a suboperad of $H_\bullet (\mathcal{D})$,
we also have a natural structure of $\mathcal{L}$-module on $H_\bullet (\mathcal{D})$.
The next proposition is a stronger version of the above lemma.
\begin{prop}
 For any $q\geqs0$, $H^\bullet (\mathcal{N}_\infty(\_,q))$ and $H_\bullet (\mathcal{D})$ are isomorphic as
  $\mathcal{L}$-modules.
\end{prop}
 %As shown below, the description of $H^\bullet (\mathcal{N}_\infty)$ is a
%consequence of Theorem \ref{cohen76}.
%In fact, we will show that $H^\bullet (\mathcal{N}_\infty(p,q))$ is isomorphic, as
%$\mathbb{Z}$-graded $k$-vector spaces, to $H_\bullet (\mathcal D(p)) = H_\bullet (C(p)) = H_\bullet (\overline{C(p)})$.
%Where $\mathcal D$ is the little discs operad, $C(p)$ is the configuration space of $p$
%points on the complex plane and $\overline{C(p)}$ denotes its real Fulton MacPherson
%compactification, see \mbox{section \ref{c(p,q)}}.
\begin{proof}
At the end of section \ref{c(p,q)} we observed that $\overline{C(p,q)}$ deformation retracts to
a stratum $S_T$ which is homeomorphic to $C(p)$.
That deformation retract takes each stratum of dimension $m$
(represented by a partially planar tree of degree $-m$) in $\overline{C(p,q)}$ to
an $m$-dimensional singular chain in $\overline{S_T} = \overline{C(p)}$.

In fact, following the notation of Theorem \ref{strata},
let $S_U$ be a stratum of dimension $m$ corresponding to a tree $U$ of degree $|U| = -m$.
Its closure $\overline{S_U}$ is a connected smooth oriented manifold with corners (topologically it is a manifold with boundary).
Let $[\overline{S_U}]$ be the relative fundamental class in
$H_m(\overline{S_U},\partial \overline{S_U})$.
For each stratum $S_U$, take a singular chain in $\overline{C(p,q)}$ representing the fundamental class
$[\overline{S_U}]$.
By composing with the contraction $\overline{C(p,q)} \rar \overline{C(p)}$,
we see that those singular chains are taken to singular chains in $\overline{C(p)}$. Hence we have a chain map:
\begin{equation}\label{Npq_chain}
 \psi_p : \mathcal{N}_\infty(p,q) \lra C_*(\overline{C(p)})
\end{equation}
and an induced map in homology
\begin{equation} \label{Npq_homology}
\Psi_p : H^\bullet (\mathcal{N}_\infty(p,q)) \lra H_\bullet (C(p)) \simeq H_\bullet (\mathcal{D}(p)),
\qquad \mbox{for each $p \geqs 1$.}\end{equation}

Since the contraction leaves the configuration of the interior points unaffected (see subsection \ref{contract}), classes representred
by trees with only spatial edges will also be unaffected. It follows that the class
$[S_{{\delta_k} \circ_i U}] = [S_{\delta_k}] \times [S_U]$ will be taken to the class $[S_{\delta_k}] \times \psi_p([S_U])$
for any spatial corolla $\delta_k \in \mathcal{L}_\infty$. Consequently,
the sequence of maps $\{ \Psi_p \}$ define a morphism of $\mathcal{L}$-modules:
\[   \Psi : H^\bullet (\mathcal{N}_\infty(\_,q)) \lra H_\bullet (\mathcal{D}).  \]

In order to show that $\Psi$ is an isomorphism, let us now construct a map from $H_\bullet (\mathcal{D}(p))$
to $H^\bullet (\mathcal{N}_\infty(p,0))$, for each $p \geqs 1$.
In case $p=1$, define the map: \[ \Phi_1 : H_\bullet (\mathcal{D}(1)) \lra H^\bullet (\mathcal{N}_\infty(1,0)) \]
by taking the identity in $e \in H_\bullet (\mathcal{D}(1))$ into
$n_{1,0} = \raisebox{-4ex}{\includegraphics{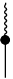}}$.

When $p=2$, the map $\Phi_2 : H_\bullet (\mathcal{D}(2)) \lra H^\bullet (\mathcal{N}_\infty(2,0))$
is defined by:
\[ \raisebox{-1.5em}{\includegraphics{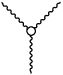}} \longmapsto
            \raisebox{-1.8em}{\includegraphics{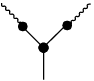}} \quad \mbox{ and }
            \raisebox{-1.5em}{\includegraphics[scale=0.8]{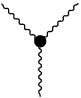}} \longmapsto
            \raisebox{-1.8em}{\includegraphics{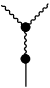}}\; . \]
Now that we have defined our maps on the operad generators of $H_\bullet (\mathcal{D})$,
we define the map
$\Phi_p : H_\bullet (\mathcal{D}(p)) \rar H^\bullet (\mathcal{N}_\infty(p,0))$, for any $p$, in the following way:
{\it \begin{enumerate}[i)]
\item if $T \in H_\bullet (\mathcal{D}(p))$ has only white vertices
(i.e., corresponds to a zero dimensional homology class), then $\Phi_p(T)$ is defined
by grafting $n_{1,0}$ to all the leaves
of the tree obtained from $T$ by making all vertices black and all edges straight;
\item extend $\Phi_p$ to the whole $H_\bullet (\mathcal{D}(p))$ so that the resulting map
\begin{equation}
\label{phi} \Phi : H_\bullet (\mathcal{D}) \rar H^\bullet (\mathcal{N}_\infty(\_,0))
\end{equation}
becomes a morphism of left modules over
$\mathcal{L} = \{ H_{n-1}(\mathcal{D}(n)) \}_{n \geqs 1}$, i.e., such that
\[ \Phi(T \circ_i l_2) = \Phi(T) \circ_i l_2, \ \mbox{for any} \ T \in H_\bullet (\mathcal{D}). \]
\end{enumerate}}
\noindent In conclusion, we have defined another morphism of $\mathcal{L}$-modules
\[ \Phi : H_\bullet (\mathcal{D}) \rar H^\bullet (\mathcal{N}(\_,0)). \]

To see that $\Phi : H_\bullet (\mathcal{D}) \rar H^\bullet (\mathcal{N}_\infty(\_,0))$ is an isomorphism, let us show
that the composition $\Psi \circ \Phi$ is the identity in $H_\bullet (\mathcal{D})$.
In fact: since both $\Psi$ and $\Phi$ are morphisms of $\mathcal{L}$-modules,
we need only to check that on generators.
Observe that $\raisebox{-2ex}{\includegraphics[scale=0.7]{n10n10.eps}}$ correspond to a zero
dimensional component of the boundary strata of $\overline{C(2,0)}$	and is taken
to the zero dimensional generator
$\raisebox{-2ex}{\includegraphics[scale=0.6]{m_2trees.eps}}
\in H_0(\mathcal{D}(2))$
under the deformation
retraction $\overline{C(2,0)} \rar \overline{C(2)}\cong \mathcal{D}(2)$ used to define $\Psi$.
On the other hand, $\raisebox{-1.8ex}{\includegraphics[scale=0.6]{l2n10.eps}}$ is homemorphic
to $S^1$ (see Figure \ref{theeye_figure} pg. \pageref{theeye_figure})
and is naturally taken to
$\raisebox{-1.8ex}{\includegraphics[scale=0.5]{l2_tree.eps}} \in H_1(\mathcal{D}(2))$
under the same deformation retraction, so $\Psi \circ \Phi = Id$. From lemma \ref{1st_isomorphism},
we know that the vector spaces $H_\bullet (\mathcal{D}(p))$ and $H^\bullet (\mathcal{N}_\infty(p,0))$ have the same
dimension for each $p \geqs 1$. It follows that $\Phi$ is in fact a bijection.

Finally we just need to observe that $H^\bullet (\mathcal{N}_\infty(\_,0))$ is naturally isomorphic as an
$\mathcal{L}$-module to $H^\bullet (\mathcal{N}_\infty(\_,q))$ for any $q \geqs 0$. The isomorphism being
induced by the grafting operation with some fixed binary planar tree $T$ with $q+1$ leaves.
\end{proof}

\begin{coro} \label{N_coro}
The cohomology $H^\bullet (\mathcal{N}_\infty)$ is
the ideal of $H^\bullet (\mathcal{OC}_\infty)$ generated by $n_{1,0} \mbox{ and } n_{0,2}.$
\end{coro}
\begin{proof}
It is immediate from the explicit definition of the $\mathcal{L}$-isomorphism
$\Phi$ that any class in $H^\bullet (\mathcal{N}_\infty(p,q))$ can be obtained
by grafting a finite number of trees of the form $n_{1,0}$ and $n_{0,2}$ followed
by grafting a finite number of the form $l_2$, i.e., by the action of
$\mathcal{L}$ on $H^\bullet (\mathcal{N}(\_,q))$.
\end{proof}

We can now prove our main result.
\begin{thm} \label{main}
The morphism of differential graded $\mathcal{L}_\infty$-modules $\mu : \mathcal{OC}_\infty \lra \mathcal{OC}$
induces an isomorphism in cohomology.
\end{thm}
\begin{proof}
It is sufficient to show that the cohomology OCHA operad $H^\bullet (\mathcal{OC}_\infty)$ and
$\mathcal{OC}$ are isomorphic as operads of graded vector spaces.
Let us first recall that the operad $\mathcal{OC}_\infty$ is decomposed as a direct sum:
$\mathcal{OC}_\infty = \mathcal{L}_\infty \oplus \mathcal{N}_\infty$,
where $\mathcal{L}_\infty$ is the operad of $L_\infty$-algebras and $\mathcal{N}_\infty$ is the is the ideal of
partially planar trees with planar root. Since the differential operator $d$ respects the direct sum
decomposition, the homology of $\mathcal{OC}_\infty$ is a direct sum:
$H^\bullet (\mathcal{OC}_\infty) = \mathcal{L}\, \oplus H^\bullet (\mathcal{N}_\infty)$.
Now we just observe that $\mathcal{L}$ is the operad generated by
$l_2$ and, from Corollary \ref{N_coro}, $H^\bullet (\mathcal{N}_\infty)$ is generated by
$n_{1,0}$ and $n_{0,2}$. The relations listed in the statement of Theorem
\ref{SC_homology} are naturally satisfied in $H^\bullet (\mathcal{OC}_\infty)$ since
they are just the homology version of the OCHA axioms.
\end{proof}

Considering that $\mathcal{OC}$
is a suboperad of $H_\bullet (\mathcal{SC})$, an interesting problem that might be pursued in a sequel
to the present paper is to extend our results to the whole operad $H_\bullet (\mathcal{SC})$. That would involve
the entire spectral sequence of $\overline{C(p,q)}$ (see also the comments at the end of \cite{Voro99}).

%\section*{Acknowledgments}
%\noindent The author wishes to thank Jim Stasheff and Murray Gerstenhaber for the kind hospitality
%during his stay as a visiting graduate student at the University of Pennsylvania (CNPq-Brasil grant
%SWE-201064/04). We are also grateful to J. Stasheff and H. Kajiura for valuable discussions.
%Some results of this paper consist of strengthened versions of some of the results present in the author's
%Ph.D. thesis \cite{Hoe06c} defended in 2006 at Unicamp under the supervision of A. Rigas
%and T. E. Barros.

%\section*{Appendix}
\appendix
\section{OCHA as a Coderivation Differential}
We say that a coderivation $\phi \in {\rm Coder}(S^c (U) \otimes T^c (V))$ is in {\bf OCHA form} if it can
be written as a sumation
\[ \phi = \sum_{n \geqs 1} \tilde g_n + \sum_{p+q \geqs 1} \tilde f_{p,q}, \]
where $\tilde g_n$ and $\tilde f_{p,q}$ denote the lifting as a coderivation of some maps:
$g_n : U^{\wedge n} \rar U$ and $f_{p,q} : U^{\wedge p} \otimes V^{\otimes q} \rar V$.
In \cite{Hoe06b} we have proven that
all coderivations in ${\rm Coder}(S^c (U) \otimes T^c (V))$ are in OCHA form for any vector spaces $U$ and $V$ over
a field $k$ of characteristic zero.

\noindent
{\bf Proposition \ref{coder_ocha}.}
{\it An OCHA structure $(L,A,\mathfrak{l},\mathfrak{n})$, in the grading and signs conventions of
defintion \ref{OCHA_def}, is equivalent to a degree one coderivation
$D \in {\rm Coder}(S^c (\downarrow \! \downarrow \!\! L) \otimes T^c ( \downarrow \!\! A))$ such that $D^2 = 0$.}
\begin{proof}
Let us begin by defining: $\tilde l_1 = - l_1$ and $\tilde n_{0,1} = - n_{0,1}$ as the differential operators
respectively on $\downarrow \! \downarrow \!\!\! L$ and on $\downarrow \!\!\! A$.
Let $D \in {\rm Coder}(S^c (\downarrow \! \downarrow \!\! L) \otimes T^c ( \downarrow \!\! A))$ be any degree one
coderivation such that $D^2 = 0$. Since any coderivation in
${\rm Coder}(S^c (\downarrow \! \downarrow \!\! L) \otimes T^c ( \downarrow \!\! A))$
is in OCHA form, $D$ is obtained by lifting maps
$\tilde l_n : (\downarrow \! \downarrow \!\!\! L)^{\otimes n} \rar \downarrow \! \downarrow \!\!\! L$
for $n \geqs 1$ and
$\tilde n_{p,q} : (\downarrow \! \downarrow \!\!\! L)^{\otimes p} \otimes (\downarrow \!\!\! A)^{\otimes q} \rar \downarrow \!\!\! A$
for $p+q \geqs 1$, where all the maps $\tilde l_n$ and $\tilde n_{p,q}$ have degree one.

Equation $D^2 = 0$ holds if and only if $\{ \tilde l_n \}_{n \geqs 1}$ satisfies the conditions of an
$L_\infty$ algebra and $\{ \tilde n_{p,q} \}_{p+q \geqs 1}$ satisfies the conditions of an OCHA as
originally defined in \cite{KS06a}:
\begin{equation}\label{ocha_shorthand}
0 = \mquad \sum_{\sigma \in \Sigma_{p+r=n}} \!\!\!
\Big( \tilde n_{1+r,m}(\tilde l_p \otimes {\bf 1}_L^{\otimes \,r} \otimes {\bf 1}_A^{\otimes \,m})
%(E(\sigma) \otimes {\bf 1}_A^{\otimes \,m}) \  + \\
 \;\; +   \sum_{\scriptstyle i+j+s=m} \!
 \tilde n_{p,i+1+j} ({\bf 1}_L^{\otimes \,p} \otimes {\bf 1}_A^{\otimes \,i} \otimes
\tilde n_{r,s} \otimes {\bf 1}_A^{\otimes \,j}) \Big) (E(\sigma) \otimes {\bf 1}_A^{\otimes \,m}).
\end{equation}
Now define maps $l_n : L^{\otimes n} \rar L$ and
$n_{p,q} : L^{\otimes p} \otimes A^{\otimes q} \rar A$, with ${\rm deg}(l_n) = 3 - 2n$ and
${\rm deg}(n_{p,q}) = 2 - 2p - q$ such that:
$\tilde l_p = \downarrow \! \downarrow l_p (\uparrow \! \uparrow)^{\otimes p} $
and $\tilde n_{p,q} = \;\downarrow n_{p,q}
(\uparrow \! \uparrow^{\otimes p} \otimes \;\uparrow^{\otimes q}).$
Thus:
\begin{multline*}
\tilde n_{1+r,m}(\tilde l_p \otimes {\bf 1}_L^{\otimes \,r} \otimes {\bf 1}_A^{\otimes \,m})
  + \!\!\!  \sum_{i+j+s=m}
\!\! \tilde n_{p,i+1+j} ({\bf 1}_L^{\otimes \,p} \otimes {\bf 1}_A^{\otimes \,i} \otimes
\tilde n_{r,s} \otimes {\bf 1}_A^{\otimes \,j}) = \\
=  \;\downarrow n_{1+r,m} (\uparrow \! \uparrow^{\otimes 1+r} \otimes \;\uparrow^{\otimes m})
(\downarrow \! \downarrow l_p (\uparrow \! \uparrow)^{\otimes p}
\otimes {\bf 1}_L^{\otimes \,r} \otimes {\bf 1}_A^{\otimes \,m}) \;+ \qquad\qquad\qquad\qquad\qquad\qquad\qquad\\ +
\sum_{i+j+s=m} \!\! \downarrow n_{p,i+1+j}
(\uparrow \! \uparrow^{\otimes p} \otimes \;\uparrow^{\otimes i+1+j})
({\bf 1}_L^{\otimes \,p} \otimes {\bf 1}_A^{\otimes \,i} \otimes
\downarrow n_{r,s} (\uparrow \! \uparrow^{\otimes r} \otimes \;\uparrow^{\otimes s})
\otimes {\bf 1}_A^{\otimes \,j}) = \\
= \;\downarrow n_{1+r,m}(l_p \otimes {\bf 1}_L^{\otimes \,r} \otimes {\bf 1}_A^{\otimes \,m})
(\uparrow \! \uparrow^{\otimes n} \otimes \;\uparrow^{\otimes m}) +
\qquad\qquad\qquad\qquad\qquad\qquad\qquad\qquad\qquad\quad\\ +
\sum_{i+j+s=m} \!\! (-1)^{js + i}\downarrow n_{p,i+1+j}
({\bf 1}_L^{\otimes \,p} \otimes {\bf 1}_A^{\otimes \,i} \otimes n_{r,s} \otimes {\bf 1}_A^{\otimes \,j})
(\uparrow \! \uparrow^{\otimes n} \otimes \;\uparrow^{\otimes m}) = \\
= \;\downarrow \Big( n_{1+r,m}(l_p \otimes {\bf 1}_L^{\otimes \,r} \otimes {\bf 1}_A^{\otimes \,m})\; + \!\!
\sum_{i+j+s=m} \!\! (-1)^{js + i}n_{p,i+1+j}
({\bf 1}_L^{\otimes \,p} \otimes {\bf 1}_A^{\otimes \,i} \otimes n_{r,s} \otimes {\bf 1}_A^{\otimes \,j})\Big)
(\uparrow \! \uparrow^{\otimes n} \otimes \;\uparrow^{\otimes m}),
\end{multline*}
where the sign $(-1)^{js + i}$ comes from the Koszul sign convention. Observing that
$\tilde l_1 = - l_1$, $\tilde n_{0,1} = - n_{0,1}$ and
$(-1)^{js + i} = (-1)^{s+i+si+ms}$, we obtain formula (\ref{geo_ocha_shorthand})
from formula (\ref{ocha_shorthand}).
\end{proof}

\section{Existence of the DG $\mathcal{L}_\infty$-module morphism}
\begin{prop}\label{L_morphism}
There is a morphism of differential graded $\mathcal{L}_\infty$-modules
$\mu : \mathcal{OC}_\infty \rar \mathcal{OC}$
extending the identity on $\mathcal{OC}$, i.e., such that the following diagram is commutative:
\begin{diagram}[h=2em]
 \mathcal{OC}_\infty &              &                  \\
 \uIntto             & \rdTTo^{\mu} &                  \\
 \mathcal{OC}        & \rTTo^{id}   &  \mathcal{OC}\ . \\
\end{diagram}
\end{prop}
\begin{proof}
 In this proof we shall omit the labels on trees because they are not crucial in the argument.

 The open-closed operad $\mathcal{OC}$ is a differential graded operad where the differential operator $\delta$
 is trivial: $\delta \equiv 0$. On the other hand, the differential operator $d$ of the OCHA operad
 $\mathcal{OC}_\infty$ is defined by formulas (\ref{dl}) and (\ref{dn}).
 We will exhibit a chain map $\mu : \mathcal{OC}_\infty \rar \mathcal{OC}$ which is also a morphism of
 $\mathcal{L}_\infty$-modules.  In other words, $\mu$ must satisfy two conditions:
 \begin{gather*}
   \mu(d T) = 0, \quad \forall T \in \mathcal{OC}_\infty \\
   \mu(l \circ_i T) = l \circ_i \mu(T), \quad \forall T \in \mathcal{OC}_\infty \ \mbox{and}\
                                                            \forall l \in \mathcal{L}_\infty.
\end{gather*}
% Consider the $\mathcal{L}$-submodule of $\mathcal{OC}_\infty$ generated by the following trees
% \[ \left \{ \raisebox{-1.5em}{\includegraphics{l2.eps}},
%             \raisebox{-1.5em}{\includegraphics{n_02.eps}},
%             \raisebox{-1.5em}{\includegraphics[scale=1.05]{n10.eps}},
%             \raisebox{-1.5em}{\includegraphics[scale=0.75]{n_11n_10.eps}} \right \}. \]
%We define $\mu$ in the following way.
\noindent
Let $\mathcal{E}$ be the $\mathcal{L}_\infty$-submodule of
$\mathcal{OC}_\infty$ generated by $\mathcal{OC}$ and by
$\raisebox{-1em}{\includegraphics[scale=0.6]{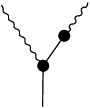}}$:
\[ \mathcal{E} =
 \left\langle \mathcal{OC}\, , \raisebox{-1em}{\includegraphics[scale=0.6]{n_11n_10.eps}} \right\rangle \]
On the generators of the submodule $\mathcal{E}$, the map $\mu$ is defined in the following way:
\[
\mu (T) = T  \quad \forall\, T \in \mathcal{OC} \quad \mbox{ and } \quad
\mu \Big( \raisebox{-1em}{\includegraphics[scale=0.6]{n_11n_10.eps}} \Big) =
- \frac{1}{2}\, \raisebox{-1.2em}{\includegraphics[scale=0.8]{l2n10.eps}}
\]
and it is extended to $\mathcal{E}$ as an $\mathcal{L}_\infty$-morphism. Finally, for any tree $T \in \mathcal{OC}_\infty$
such that $T \notin \mathcal{E}$, we define $\mu(T) = 0$. We thus have an $\mathcal{L}_\infty$-morphism:
\[ \mu : \mathcal{OC}_\infty \rar \mathcal{OC}. \]
It remains to show that $\mu$ is a chain map, i.e., that $\mu(d\, T) = 0$ for any tree $T \in \mathcal{OC}_\infty$.
Given any tree $T \in \mathcal{OC}_\infty$, $d\, T$ is a summation of trees. By the definition of $\mu$,
if $T$ is such that $d\, T$ has no components in $\mathcal{E}$, then $\mu(d\, T) = 0$.
Hence, we just need to consider those elements $T \in \mathcal{OC}_\infty$ such that
$d\, T$ has some component in $\mathcal{E}$. Such elements form an $\mathcal{L}_\infty$-submodule of
$\mathcal{OC}_\infty$ which will be denoted by $\mathcal{E}'$. More precisely:
\[ \mathcal{E}' := \{ T \in \mathcal{OC}_\infty \,:\, d\, T = T_1 + T_2, \quad T_1 \in \mathcal{E}, T_1 \neq 0 \}. \]

Any tree $T$ is obtained by grafting a finite number of corollae which we call the {\it irreducible components}
of $T$.
Recall that, for $n \geqs 3$, the $\mathcal{L}_\infty$-module action of $l_n \in \mathcal{L}_\infty$
on any element of $\mathcal{OC}$
is zero since that action is defined through the quasi-ismorphism $\mu : \mathcal{L}_\infty \rar \mathcal{L}$,
and $\mu(l_n) = 0$ for $n \geqs 3$.
From the definition of $\mu : \mathcal{OC}_\infty \rar \mathcal{OC}$ and the definition of the
$\mathcal{L}_\infty$-module structure on $\mathcal{OC}$, one can see that the irreducible components of any tree
$T \in \mathcal{E}'$ such that $\mu(d T) \neq 0$ could only be one of the following corollae:
\[ \left \{  \raisebox{-1.5em}{\includegraphics[scale=0.95]{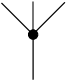}},
             \raisebox{-1.5em}{\includegraphics{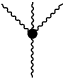}},
             \raisebox{-1.5em}{\includegraphics{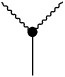}},
             \raisebox{-1.5em}{\includegraphics{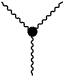}},
             \raisebox{-1.5em}{\includegraphics[scale=0.95]{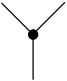}},
             \raisebox{-1.5em}{\includegraphics{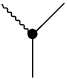}},
             \raisebox{-1.5em}{\includegraphics{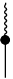}} \right \}. \]
Consequently, we just need to check that $\mu (d\, T) = 0$ where $T$ is any of the above corollae.
In the case of $T = \raisebox{-.8em}{\includegraphics[scale=0.7]{n_20.eps}}: \quad$
$
      \mu \Big( d\, \raisebox{-1.5em}{\includegraphics[scale=0.9]{n_20.eps}} \Big)
    = \mu \Big( \;\raisebox{-1.5em}{\includegraphics[scale=0.7]{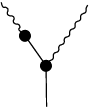}}
      + \raisebox{-1.5em}{\includegraphics[scale=0.7]{n_11n_10.eps}}
      + \raisebox{-1.5em}{\includegraphics[scale=0.87]{l2n10.eps}} \;\Big)
    = - \raisebox{-1.5em}{\includegraphics[scale=0.87]{l2n10.eps}}
      + \raisebox{-1.5em}{\includegraphics[scale=0.87]{l2n10.eps}} = 0,
$
since by definition we have:
 $\mu (\raisebox{-1.2em}{\includegraphics[scale=0.6]{n_11n_10.eps}}) =
  - \frac{1}{2}\, \raisebox{-1.2em}{\includegraphics[scale=0.8]{l2n10.eps}}$,
 $\quad \mu (\raisebox{-1.2em}{\includegraphics[scale=0.8]{l2n10.eps}}) =
       \raisebox{-1.2em}{\includegraphics[scale=0.8]{l2n10.eps}}$
and because the wiggly edges are spatial, we also have:
$\raisebox{-1.2em}{\includegraphics[scale=0.6]{fn_11n_10.eps}} = \raisebox{-1.2em}{\includegraphics[scale=0.6]{n_11n_10.eps}}$.
The other corollae can be handled similarly. \qedhere
\end{proof}

\begin{figure}[cp]
\centering
 {\includegraphics[scale=0.5]{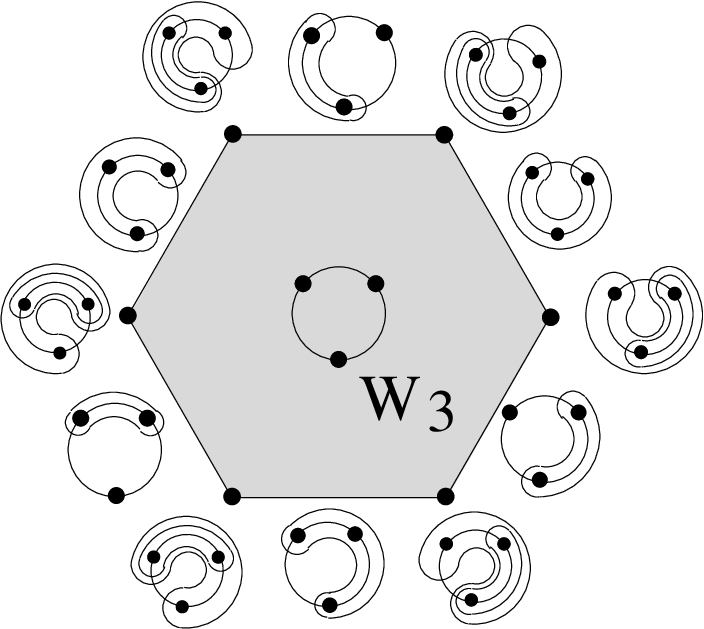}} \hspace{4em}
 {\includegraphics[scale=0.5]{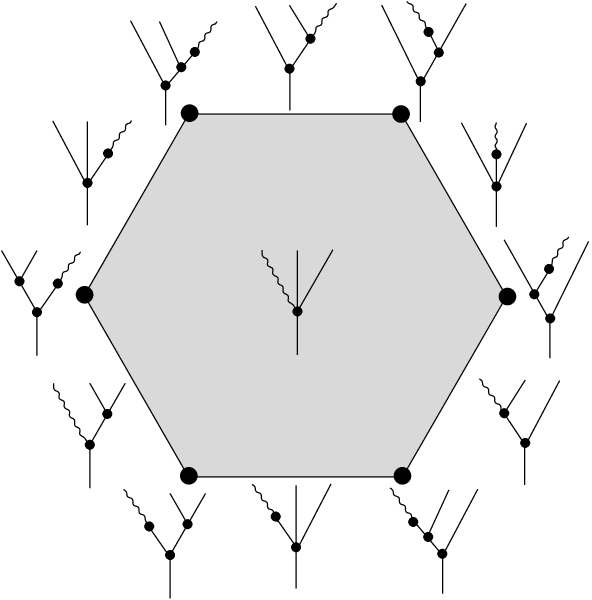}}
\caption{Cyclohedron $\overline{C(1,2)}$ and its cells labelled by circular
bracketings and by trees.}
\label{w3}
\vspace*{2em}
\end{figure}

\begin{figure}[cp]
\centering
\raisebox{-10em}{\input{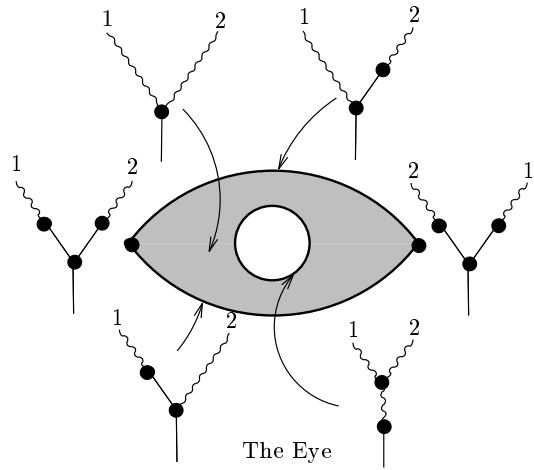}}
\caption{The space $\overline{C(2,0)}$ = ``The Eye'' and its boundary
strata labelled by trees.}
\label{theeye_figure}
\vspace*{2em}
\end{figure}

\begin{figure}[cp]
\centering
\input{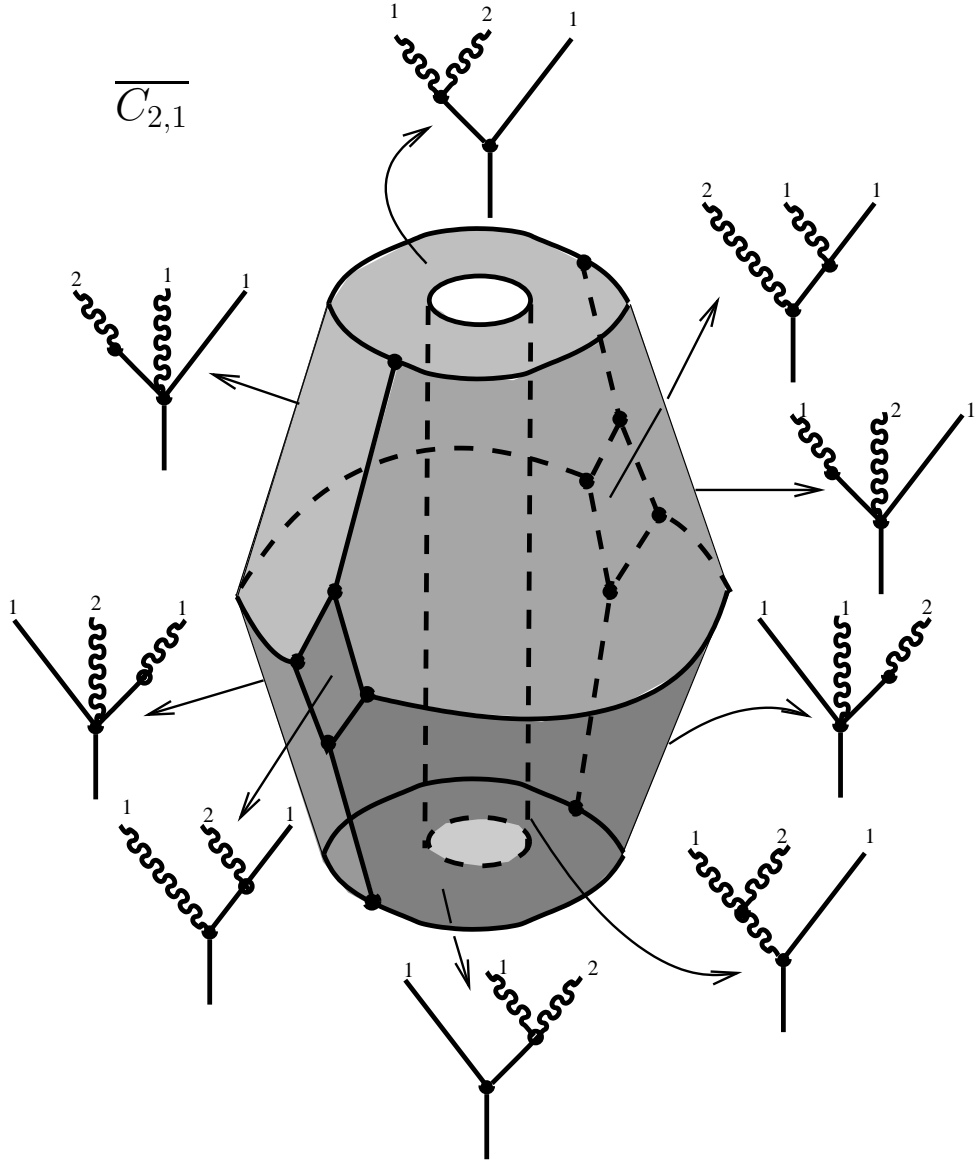}
\caption{ The space $\overline{C(2,1)}$, which is topologically
equivalent to a solid torus, and its codimension 1 boundary components labelled by partially planar trees.}
\label{c21}
\end{figure}

\bibliography{eduII}

\end{document}